\documentclass[reqno,10pt]{amsart}

\usepackage{a4wide}
\usepackage{color}
\usepackage[T2A]{fontenc}
\usepackage[utf8]{inputenc}
\usepackage[english]{babel}
\numberwithin{equation}{section}
\usepackage[colorlinks,citecolor=blue,linkcolor=red]{hyperref}
\usepackage[noabbrev,capitalise,nameinlink]{cleveref}
\hypersetup{colorlinks={true},linkcolor={red},citecolor=blue}
\usepackage{csquotes}

\makeatletter
\newcommand{\biggg}{\bBigg@{3}}
\newcommand{\Biggg}{\bBigg@{3.5}}
\newcommand{\bigggg}{\bBigg@{4}}
\newcommand{\bdot}{\text{$\,\cdot\,$}}
\newcommand{\lims}{\varlimsup}
\DeclareMathOperator*{\esup}{ess\,sup}
\DeclareMathOperator*{\einf}{ess\,inf}
\newcommand{\limi}{\varliminf}
\DeclareMathOperator*{\supp}{supp}

\newcommand{\D}{\mathrm{d}}
\newcommand{\Leb}{\mathrm{L}}

\newcommand{\dist}{\mathrm{dist}}

\makeatother

\usepackage{amsmath, amsfonts, amssymb, amsthm, mathtools, bbm, esint, mathrsfs, nccmath}
\usepackage[backend=biber,style=numeric-comp,sortcites,giveninits,doi=false,isbn=false,url=false,eprint=false]{biblatex}

\usepackage{resmes}

\AtBeginBibliography{\footnotesize\sloppy}

\DeclareNameAlias{author}{family-given}

\theoremstyle{plain}
\newtheorem{theo}{Theorem}[section]
\newtheorem{lem}[theo]{Lemma}
\newtheorem{prop}[theo]{Proposition}
\newtheorem{cor}{Corollary}

\theoremstyle{definition}
\newtheorem{defi}{Definition}[section]

\theoremstyle{remark}
\newenvironment{rema} {\begin{proof}[Remark]} {\end{proof}}

\defbibheading{bibliography}[\refname]{%
  \section*{#1}
}

\usepackage{relsize}

\title[Asymptotic relations of the BBM type for mappings between singular spaces]{Asymptotic relations\\
of the Bourgain-Brezis-Mironescu type\\
for mappings\\
between singular spaces}

\author[R. D. Oleinik]{Roman D. Oleinik}
\address{Moscow Institute of Physics and
Technology, 9 Institutskiy per., Dolgoprudny, Moscow Region, 141701, Russian Federation}
\email{oleinik.r@phystech.edu}

\begin{document}
\allowdisplaybreaks

\begin{abstract}
We explore the asymptotic behavior of families of Bourgain-Brezis-Mironescu type nonlocal functionals for mappings from metric measure spaces to arbitrary metric spaces. As the first outcome, we obtain a characterization of Sobolev maps and of maps of bounded variation via such functionals. As the second outcome, we establish precise expressions of the limits of such functionals for Sobolev maps. All this provides an extension of several Bourgain-Brezis-Mironescu type results to the entirely singular setting.
\end{abstract}

\maketitle
\tableofcontents

\section{Introduction}
\subsection{Overview}
Investigating in \cite{BBM01,BBM02} the limiting behavior of a special class of nonlocal functionals on Euclidean domains, J. Bourgain, H. Brezis, and P. Mironescu established a new remarkable attribute of Sobolev functions. It was shown that, given a smooth bounded domain $\Omega\subseteq \mathbb{R}^d$ with $d\in \mathbb{N}$ and a certain family $(\varrho_{\delta})_{\delta\in (0,1)}$ of mollifiers on $\mathbb{R}^d$, understood as functions $\mathbb{R}^d\times \mathbb{R}^d\to [0,+\infty)$, the following formula holds for any $u\in \Leb^p(\Omega)$ with $p\in (1,+\infty)$:
\begin{equation}
\lim\limits_{\delta\searrow 0}\int\limits_{\Omega\times \Omega} \frac{\big|u(x')-u(x)\big|^p}{\|x'-x\|^{p}}\varrho_{\delta}(x,x') \D (\mathcal{L} \otimes\mathcal{L})(x,x')=C(p,d)\mathrm{E}_{p}[u](\Omega),
\label{eq: Seminorm}
\end{equation}
where $\|\bdot\|$ and $\mathcal{L}$ are the Euclidean norm and the Lebesgue measure, respectively, in $\mathbb{R}^d$, while $C(p,d)$ is a positive constant depending only on $p$ and $d$. In the above expression, $\mathrm{E}_{p}[u](\Omega)$ stand for the $p$-energy of $u$ on $\Omega$ defined by
\begin{equation*}
\mathrm{E}_{p}[u](\Omega)\coloneqq \begin{dcases}\int\limits_{\Omega} \big\|(\nabla u)(x)\big\|^p \D \mathcal{L}(x),\quad u\in \mathrm{W}^{1,p}_{Loc}(\Omega),\\
+\infty, \quad \text{otherwise},
\end{dcases}
\end{equation*}
where $\nabla u$ is the weak gradient of $u$. An analogous property was then discovered, by J. D{\'a}vila in \cite{D02}, for functions of bounded variation. Namely, it was obtained that \eqref{eq: Seminorm}, with $p=1$, holds for any $u\in \Leb^1(\Omega)$ if one assigns $\mathrm{E}_{1}[u](\Omega)$, the $1$-energy of $u$ on $\Omega$, to be the, possibly infinite, total variation of $u$ on $\Omega$. Subsequently, these results were revisited in the Euclidean setting several times. We refer to \cite{PSV19, P04, L14, L14.2, PSV19, NPSV20, NS19, BN16, BN16.2, BN20, LS11, LS14, BVSY21, N11, BSVSY22, GT22, GT20, CS19, BCCS22, CS22, BSVSY21} as a non-exhaustive list of papers dealing with related issues. Results of this kind are now commonly labeled as ones of BBM-type after Bourgain-Brezis-Mironescu.

In recent decades, studies on first-order Sobolev calculus in the singular setting of metric measure spaces received particular importance and prevalence. A metric measure space is usually meant as a triple $(\mathrm{X}, \mathsf{d}, \mathfrak{m})$ consisting of a separable metric space $(\mathrm{X}, \mathsf{d})$ endowed with a locally finite Borel-regular outer measure $\mathfrak{m}$. For great sources on the topic, we send the reader to the monographs \cite{BB11} by A. Bj{\"o}rn and J. Bj{\"o}rn, \cite{HKST15} by J. Heinonen, P. Koskela, N. Shanmugalingam, and J. T. Tyson, and \cite{GP20} by N. Gigli and E. Pasqualetto. In connection with aforesaid, it is natural to wonder whether BBM-type results can be obtained in the described situation, at least under some reasonable conditions on the underlying space and on the family of mollifiers.

Since no additional regularity structure on metric measure spaces is a priori assumed, one can hope in general to achieve only an equivalence, up to some universal structural constants, between, both appropriately extended to the new setting, the left-hand side, with the limit replaced by the lower and the upper limits, and the right-hand side of \eqref{eq: Seminorm}. The first result of this type was provided in \cite{M15} by V. Munnier. It was, in turn, revisited by S. Di Marino and M. Squassina in \cite{DS19}. The latest result is obtained in \cite{LPZ22} by P. Lahti, A. Pinamonti, and X. Zhoue. It was considered there a complete connected metric measure space $(\mathrm{X}, \mathsf{d}, \mathfrak{m})$ supporting doubling and $(p,p)$-Poincar{\'e} inequalities for some $p\in [1,+\infty)$, together with a family $(\varrho_{\delta})_{\delta\in (0,1)}$ of mollifiers satisfying certain, fairly flexible, properties. It was established that if $\Omega\subseteq \mathrm{X}$ is a strong $p$-extension domain, then for any $u\in \Leb^p(\Omega)$ one has
\begin{gather}
\frac{1}{C}\mathrm{E}_p[u] (\Omega)\leq  \limi\limits_{\delta\searrow 0}\int\limits_{\Omega\times \Omega} \frac{\big|u(x')-u(x)\big|^p}{\big(\mathsf{d}(x,x')\big)^{p}} \varrho_{\delta}(x,x') \D(\mathfrak{m}\otimes \mathfrak{m})(x,x')\leq\notag\\
\leq \lims\limits_{\delta\searrow 0}\int\limits_{\Omega\times \Omega} \frac{\big|u(x')-u(x)\big|^p}{\big(\mathsf{d}(x,x')\big)^{p}} \varrho_{\delta}(x,x') \D(\mathfrak{m}\otimes \mathfrak{m})(x,x')\leq C \mathrm{E}_p[u]( \Omega),
\label{eq: DSRes}
\end{gather}
where $C$ is a universal positive constant depending only on $p$, on the doubling constant and the Poincar{\'e} constants of the underlying space, and on the constant associated with the mollifiers.

In order to obtain a relation that generalizes \eqref{eq: Seminorm} to the singular context in full measure, one needs to consider such a metric measure space that "has a nice local structure" in some appropriate sense. The earliest result of this sort was proposed by W. G{\'o}rny in \cite{G22}. It was then extended by B.-X. Han and A. Pinamonti in \cite{HP21}. A regularity condition for a complete connected metric measure space $(\mathrm{X}, \mathsf{d}, \mathfrak{m})$ that was adopted there is that $\mathfrak{m}$-a.e. $x\in \mathrm{X}$ has a well-behaved tangent cone in the pointed measured Gromov-Hausdorff sense. Under the additional assumption that the source space supports doubling and $(1,p)$-Poincar{\'e} inequalities for some $p\in [1,+\infty)$, it was established that whenever $u\colon \mathrm{X}\to \mathbb{R}$ is a compactly supported Lipschitz function and $(\varrho_{\delta})_{\delta\in (0,1)}$ is a family of mollifiers of a special type, then
\begin{gather}
\lim\limits_{\delta\searrow 0} \int\limits_{\mathrm{X}\times \mathrm{X}}\frac{\big|u(x')-u(x)\big|^p}{\big(\mathsf{d}(x,x')\big)^{p}}\varrho_{\delta}(x,x')\D(\mathfrak{m}\otimes \mathfrak{m})(x,x')
=\int\limits_{\mathrm{X}}\epsilon_p[u](x) \D \mathfrak{m}(x),
\label{eq: HPRes}
\end{gather}
where $\epsilon_p[u]$ is a non-negative function constructed through a specific differential of $u$ defined on the corresponding tangent cones. Additionally, it was shown, under some auxiliary mild conditions, that the above relation holds for general Sobolev functions as well, but only in the case $p>1$ and only in the limit sense with respect to approximations by Lipschitz functions.

As one can notice, all the mentioned results concern functions, i.e. real-valued maps. At the same time, Sobolev theory of metric-valued maps, i.e. maps that take values in an arbitrary metric space, is extremely significant for modern geometric analysis. Examples of articles on the topic are the pioneer works \cite{A90} by L. Ambrosio, \cite{KS93} by N. J. Korevaar and R. M. Schoen, and \cite{R97} by Y. G. Reshetnyak, as well as the recent studies \cite{GPS20} by N Gigli, E. Pasqualetto, and E. Soultanis, \cite{GT21} by N. Gigli and A. I. Tyulenev, and \cite{MS22} by A. Mondino and D. Semola. Therefore, it is of a particular interest to extend the relations presented in \eqref{eq: DSRes} and \eqref{eq: HPRes} to the outlined setting as well. In this regard, the goal of the current paper is to provide possible generalizations of the results by Lahti-Pinamonti-Zhoue and by Han-Pinamonti to the metric-valued case. And although such an extension for the former can hypothetically be derived somehow from the corresponding real-valued result, for the latter, in turn, it is hardly feasible, since the real-valued statement in this case relies heavily on a Rademacher-type theorem that is known to be not valid for metric-valued maps in general. So in order to encompass the whole picture, we investigate both issues from the ground up within our manuscript.
\subsection{Results} Let us formulate now the main subject of study in the article more specifically. Given a metric measure space $(\mathrm{X}, \mathsf{d}, \mathfrak{m})$, a parameter $p\in [1,+\infty)$, a family $(\rho_{\delta})_{\delta\in (0,1)}$ of Borel functions $\mathrm{X}\times \mathrm{X}\to [0,+\infty)$, an open set $\Omega\subseteq \mathrm{X}$, a metric space $(\mathrm{Y}, \mathsf{d}_{\mathrm{Y}})$, and an essentially separably valued Borel map $f\colon \Omega\to \mathrm{Y}$, we aim to analyze the behavior as $\delta\searrow 0$ of the functions
\begin{gather}
 \int\limits_{\Omega} \Big(\mathsf{d}_{\mathrm{Y}}\big(f(\bdot),f(x')\big)\Big)^p\rho_{\delta}(\bdot,x')\D \mathfrak{m}(x'), \quad\text{$\delta\in (0,1)$},
\label{eq: IntSub}
\end{gather}
under some reasonable additional conditions on all the objects involved. Note that in the expression above we  use, and will keep doing this throughout the paper because of convenience reasons, a slightly different, as opposed to conventional, representation for the integrand. Thus, it is the members of families such as $(\rho_{\delta})_{\delta\in (0,1)}$ that will hereafter be referred to as mollifiers, provided that they enjoy some satisfactory properties. To accomplish the stated purposes, we restrict ourselves only to those metric measure spaces that are also locally complete and connected. Besides this, some standard geometric conditions, namely doubling and Poincar{\'e} ones, need to be imposed on the source space. And since all reasoning in our context will be local to some extent, we shall deal with their appropriate relaxed variations, close to the weakest possible under which one can expect to reach the stated goals.

In our first result, we formulate a sufficient condition for $\mathrm{E}_p[f](\Omega)$, the $p$-energy of $f$ on $\Omega$, to be equivalent, up to structural multiplicative constants, to the lower and the upper limits as $\delta\searrow 0$ of the integrals over $\Omega$ of the functions in \eqref{eq: IntSub}. As an outcome, we obtain a metric-valued analog to the result by Lahti-Pinamonti-Zhoue. The precise assertion is as follows.
\begin{theo}\label{theo: Res1}
Let $(\mathrm{X}, \mathsf{d},\mathfrak{m})$ be a locally complete connected metric measure space. Let $p\in [1,+\infty)$, suppose $(\mathrm{X}, \mathsf{d},\mathfrak{m})$ is internally doubling and supports an internal $p$-Poincar{\'e} inequality, let $(\rho_{\delta})_{\delta\in (0,1)}$ be a $p$-admissible family. Then there exists a constant $C\in (0,+\infty)$ depending only on $p,C_D,C_P, \lambda_P,C_M$ such that the following holds. Let $\Omega\Subset \mathrm{X}$ be an open set having the strong $p$-extension property, let $(\mathrm{Y},\mathsf{d}_{\mathrm{Y}})$ be a metric space, let $f\in \Leb^p(\Omega, \mathrm{Y})$. The estimate below holds:
\begin{gather}
\frac{1}{C}\mathrm{E}_p[f](\Omega)\leq\limi\limits_{\delta\searrow 0} \int\limits_{\Omega\times \Omega} \Big(\mathsf{d}_{\mathrm{Y}}\big(f(x),f(x')\big)\Big)^p\rho_{\delta}(x,x') \D (\mathfrak{m}\otimes \mathfrak{m})(x,x') 
\leq \notag\\
\leq\lims\limits_{\delta\searrow 0} \int\limits_{\Omega\times \Omega} \Big(\mathsf{d}_{\mathrm{Y}}\big(f(x),f(x')\big)\Big)^p\rho_{\delta}(x,x') \D (\mathfrak{m}\otimes \mathfrak{m})(x,x')  \leq C\mathrm{E}_p[f](\Omega).
\label{eq: Res1}
\end{gather}
\end{theo}
\noindent The above statement, when applied to the real-valued case, essentially turns into the result by Lahti-Pinamonti-Zhoue. The only contrast, apart from that we use relaxed variants of doubling and Poincar{\'e} conditions, is that we consider different, but still quite general, assumptions on mollifiers. To prove the theorem, we follow ideas not too far from those exploited within the proof of the corresponding real-valued statement. Nevertheless, our proof includes, in addition to more involved calculations, some procedures specific to working with metric-valued maps, like those related to energies of such maps. As an intermediary result, we will also establish an equivalence similar to the above one, but which is fulfilled without integrability assumptions on the map and without extension assumptions on the open set, which may turn out to be useful in some situations.

Within the second result, we demonstrate a situation when the limit as $\delta\searrow 0$ of the functions from \eqref{eq: IntSub} actually exists, both in the integral sense and almost everywhere. This statement extends, but only partially, that of Han-Pinamonti to metric-valued maps. The explicit formulation is as follows.
\begin{theo}\label{theo: Res2}
Let $(\mathrm{X}, \mathsf{d}, \mathfrak{m})$ be a locally complete connected metric measure space. Suppose $(\mathrm{X}, \mathsf{d}, \mathfrak{m})$ is strongly rectifiable. Let $p\in [1,+\infty)$, suppose $(\mathrm{X}, \mathsf{d}, \mathfrak{m})$ is internally doubling and supports an internal $p$-Poincar{\'e} inequality, let $(\rho_{\delta})_{\delta\in (0,1)}$ be a strongly $p$-admissible family. Let $\Omega\Subset\mathrm{X}$ be an open set, let $(\mathrm{Y},\mathsf{d}_{\mathrm{Y}})$ be a metric space, let $f\in \mathrm{W}^{1,p}(\mathrm{X},\mathrm{Y})$. Then there exists a non-negative function $\epsilon_p[f]\in \Leb^1(\Omega)$ satisfying the properties below:
\begin{gather}
\lim\limits_{\delta\searrow 0}\int\limits_{ \Omega}\biggg|\int\limits_{\Omega} \Big(\mathsf{d}_{\mathrm{Y}}\big(f(x),f(x')\big)\Big)^p\rho_{\delta}(x,x') \D\mathfrak{m}(x')-\epsilon_p[f](x)\biggg|\D \mathfrak{m}(x)=0,\;\label{eq: Res2}\\
\lim\limits_{\delta\searrow 0}\biggg|\int\limits_{\Omega} \Big(\mathsf{d}_{\mathrm{Y}}\big(f(x),f(x')\big)\Big)^p\rho_{\delta}(x,x') \D\mathfrak{m}(x')-\epsilon_p[f](x)\biggg|=0\quad \text{for $\mathfrak{m}$-a.e. $x\in \Omega$}.\label{eq: Res2'}
\end{gather}
\end{theo}
\noindent Once again, besides it is related exactly to metric-valued maps, our statement differs from that of Han-Pinamonti in requirements for mollifiers, as well as in doubling and Poincar{\'e} conditions. A more relevant distinction pertains the choice of a regularity property for the underlying space. The one we deal with is the strong rectifiability property introduced by N. Gigli and E. Pasqualetto in \cite{GP22}. Spaces satisfying it are worthy of interest because they encompass the $\mathrm{RCD}$-space, which, in turn, are highly relevant for modern geometric measure theory. For a thorough discussion on the topic, we send the reader to \cite{A18} by L. Ambrosio. At the same time, the corresponding class of strongly rectifiable spaces does not cover all of those ones that fall into the scope of the result by Han-Pinamonti. And the reason we focused specifically on such spaces is that for Sobolev metric-valued maps defined on them there exist so-called energy densities. This powerful result, established by N. Gigli and A. I. Tyulenev in \cite{GT21}, is a cornerstone component for the proof of our theorem. In particular, it is this result that allows the relations from the theorem to be derived, even for Sobolev maps, directly, rather than through approximating by Lipschitz maps. But it, however, is established in the original paper only for the case $p>1$, so, as an intermediate part of the proof, we extend it additionally to the case $p=1$. Due to this, we become able to cover the whole range of exponents within our theorem, even though this will necessitate the use of slightly more involved techniques. And finally, in connection with everything that was just said, it is worth emphasizing the following. One collateral outcome of our proof is that under the assumption that the energy density of the map exists, the conclusion of the stated theorem is fulfilled even without the requirement for the source space to be strongly rectifiable. Thereby, our result is actually applicable to more general situations.

The paper is organized as follows. We give all necessary definitions in \cref{ss: Prelim}. After this, in \cref{ss: Proof}, we provide the proofs for our main theorems. And lastly, in \cref{ss: App}, we highlight some concrete consequences of our results.

{\bf Acknowledgments.} The author is indescribably grateful to his Scientific Advisor Alexander I. Tyulenev for inspiring to write this article, for useful tips on the topic, and for helping to correct mistakes and typos in the paper.
\section{Preliminaries}\label{ss: Prelim}
\subsection{Basic notation}\label{ss: GenNot} We start by introducing some notation that we will exploit later on.

Until the end of the manuscript, the following data is assumed fixed:
\begin{itemize}
\item[$\bullet$] a separable locally complete connected metric space $(\mathrm{X}, \mathsf{d})$;
\item[$\bullet$] a locally finite Borel-regular outer measure $\mathfrak{m}$ on $(\mathrm{X}, \mathsf{d})$.
\end{itemize}
The metric measure space $(\mathrm{X} ,\mathsf{d}, \mathfrak{m})$ will be denoted just by $\mathrm{X}$ in the sequel. To avoid discussing separately degenerate cases, we initially assume that our space is of positive diameter and of positive measure.

In the case when $\mathrm{X}$ is not complete, it will be convenient to consider also its metric completion, which we denote by $\overline{\mathrm{X}}$. Since $\mathrm{X}$ is locally complete, we may assume without loss of generality that it is embedded into $\overline{\mathrm{X}}$ as a dense open subset.

Everywhere below, given $\alpha\in [0,+\infty)$, we use the notation $\mathbb{R}_{\geq \alpha}\coloneqq [\alpha,+\infty)$ and $\mathbb{R}_{>\alpha}\coloneqq (\alpha,+\infty)$. Given $d\in \mathbb{N}$, by $\mathcal{L}^d$ and $B^d$ we denote the Lebesgue measure and the standard unit ball, respectively, in $\mathbb{R}^d$.

Given $S_1,S_2\subseteq \overline{\mathrm{X}}$, let $\dist(S_1,S_2)$ be the distance between $S_1$ and $S_2$ with respect to the metric on $\overline{\mathrm{X}}$. Given $x\in \mathrm{X}$, $r\in \mathbb{R}_{> 0}$, and $S\subseteq \mathrm{X}$, we put
\begin{gather*}
B(x,r)\coloneqq \big\{x'\in \mathrm{X} \mid \mathsf{d}(x,x')\leq r\big\}, \quad B(S,r)\coloneqq \Big\{x'\in \mathrm{X} \bigm| \dist\big(S,\{x'\}\big)\leq r\Big\}.
\end{gather*}

Given $S_0,S\subseteq \mathrm{X}$ with $S\subseteq S_0$, we write $S\Subset S_0$ whenever $\dist\big(S,\overline{\mathrm{X}}\backslash S_0\big)>0$, in which case the set $S$ is said to be a well-contained subset of $S_0$.

\begin{rema}
The reason we need the just introduced notion is as follows. Given an open subset of our space, we will occasionally want to be able to find its "nondegenerate" neighborhood on which our measure is at least locally finite. As soon as the space is complete, it is not a problem. Furthermore, every open subset of the space is its well-contained subset in this case. But once it is incomplete, such neighborhoods, even if understood with respect to the completion of the space, may not exist because of the unsatisfactory behavior of the measure. But by requiring the open set to be a well-contained subset of the space, we exclude such a situation from our consideration. Thus, the use of such a notion allows us to describe both cases simultaneously.
\end{rema}

Given $S\subseteq \mathrm{X}$, by $\chi_{S}$ we mean the characteristic function of $S$. Given arbitrary sets  $S_1,S_2,S$ with $S\subseteq S_1$ and a map $f\colon S_1\to S_2$, by $f\big|_{S}$ we denote the restriction of $f$ to $S$. Given arbitrary sets $S_0,S$ with $S\subseteq S_0$ and an outer measure $\mu$ on $S_0$, let $\mu \resmes_{S}$ stand for the restriction of $\mu$ to the family of subsets of $S$.

Given a metric space $(\mathrm{Y}, \mathsf{d}_{\mathrm{Y}})$, let $\mathrm{BLip}_1(\mathrm{Y})$ denote the set of all bounded functions $\mathrm{Y}\to \mathbb{R}$ that are also $1$-Lipschitz.

From here on, given a set $S\subseteq \mathrm{X}$, a metric space $(\mathrm{Y},\mathsf{d}_{\mathrm{Y}})$, and a map $f\colon S\to \mathrm{Y}$, by $\mathsf{d}_f$ we mean the function $S\times S\to \mathbb{R}$ defined as
\begin{equation*}
\mathsf{d}_f(x,x')\coloneqq \mathsf{d}_{\mathrm{Y}}\big(f(x),f(x')\big).
\end{equation*}

Let $p\in \mathbb{R}_{\geq 1}$, let $E\subseteq \mathrm{X}$ be a Borel set, let $(\mathrm{Y},\mathsf{d}_{\mathrm{Y}})$ be a metric space. Given $S\subseteq \mathrm{Y}$, let $\mathfrak{B}(E,S)$ be the collection of all, equivalence classes up to $\mathfrak{m}$-a.e. equality of, essentially separably valued Borel maps $E\to \mathrm{Y}$ whose image is contained in $S$. The same notation, but understood with respect to $(\mathfrak{m}\otimes \mathfrak{m})$, will be used for Borel subsets of $\mathrm{X}\times \mathrm{X}$. By $\Leb^p(E)$ we mean the space of all functions belonging to $\mathfrak{B}(E,\mathbb{R})$ that are $p$-integrable on $E$. Next, let $\Leb^p_{Loc}(E)$ denote the space of all functions $u\in \mathfrak{B}(E, \mathbb{R})$ with the property that for any $x\in E$ there is $r\in \mathbb{R}_{>0}$ such that $u$ is $p$-integrable on $\big(E\cap B(x,r)\big)$. Finally, let $\Leb^p(E,\mathrm{Y})$, and respectively $\Leb^p_{Loc}(E,\mathrm{Y})$, stand for the space of all maps $f\in \mathfrak{B}(E,\mathrm{Y})$ for which there is $y\in \mathrm{Y}$ such that $\mathsf{d}_{\mathrm{Y}}\big(f(\bdot),y\big)\in \Leb^p(E)$, and respectively $\mathsf{d}_{\mathrm{Y}}\big(f(\bdot),y\big)\in \Leb^p_{Loc}(E)$.

The following simple statement provides a property necessary for the further discussion which is enjoyed by essentially separably valued Borel maps.
\begin{prop}\label{prop: LipSup}
Let $E\subseteq \mathrm{X}$ be a Borel set, let $(\mathrm{Y},\mathsf{d}_{\mathrm{Y}})$ be a metric space, let $f\in \mathfrak{B}(E, \mathrm{Y})$. Then there exist a family $(\phi_j)_{j\in \mathbb{N}}\subseteq \mathrm{BLip}_1(\mathrm{Y})$ and a set $Q_{SV}\subseteq E$ with $\mathfrak{m}(E\backslash Q_{SV})=0$ such that
\begin{equation}
\mathsf{d}_f(x,x')=\sup\limits_{j\in \mathbb{N}}\mathsf{d}_{\phi_j\circ f}(x,x')\quad \text{for all $x,x'\in Q_{SV}$}.
\label{eq: LipSup}
\end{equation}
In particular, we have that $\mathsf{d}_f\in \mathfrak{B}\big(E\times E, \mathbb{R}_{\geq 0}\big)$.
\begin{proof}
We give an outline of the proof only. One needs to take the distance functions to points from a countable dense subset of the essential image of $f$ and consider then their truncations. This immediately gives the first part. The second part, in turn, follows from the first one via standard measure-theoretic arguments.
\end{proof}
\end{prop}
\begin{rema}
In fact, all the subsequent reasoning shows that one can take as the core class of maps all weakly Borel ones that satisfy the property stated in the proposition above, and not only essentially separably valued Borel ones. This would not lead however to any relevant generalization, as one can deduce from what follows.
\end{rema}

In the sequel, we will apply, both for series and integrals, the classical convergence theorems, namely, Fatou's lemma, reverse Fatou's lemma, the monotone convergence theorem, and the dominated convergence theorem. Besides this, we will use the generalized version of the latter, which we formulate here in the following way.
\begin{prop}\label{prop: GenLebTh}
Let $E\subseteq \mathrm{X}$ be a Borel set, let $\Gamma\subseteq \mathbb{R}$ be a set having zero as a limit point, let $(u_{\gamma})_{\gamma\in \Gamma},(v_{\gamma})_{\gamma\in \Gamma}\subseteq \mathfrak{B}(E,\mathbb{R}_{\geq 0})$. Suppose the following conditions are satisfied:
\begin{itemize}
\item[$\bullet $] for each $\gamma\in \Gamma$ one has that $u_{\gamma}(x)\leq v_{\gamma}(x)$ for $\mathfrak{m}$-a.e. $x\in E$;
\item[$\bullet $] it holds that $(u_{\gamma})_{\gamma\in \Gamma}$ converges $\mathfrak{m}$-a.e. on $E$ to some function $u\in \mathfrak{B}(E,\mathbb{R}_{\geq 0})$ as $\gamma \rightarrow 0$;
\item[$\bullet $] it holds that $(v_{\gamma})_{\gamma\in \Gamma}$ converges in $\Leb^1(E)$ to some function $v\in \Leb^1(E)$ as $\gamma \rightarrow 0$.
\end{itemize}
Then $u\in\Leb^1(E)$ and $(u_{\gamma})_{{\gamma}\in \Gamma}$ converges in $\Leb^1(E)$ to $u$ as $\gamma \rightarrow 0$.
\begin{proof}
We refer to \cite[Lemma 3.14]{GT21} for the proof.
\end{proof}
\end{prop}

Let $\Omega\subseteq \mathrm{X}$ be an open set. By $\mathrm{Lip}_{Loc}(\Omega)$ we mean the set of all functions $u\colon \Omega \to \mathbb{R}$ with the property that for any $x\in \Omega$ there is $r\in \mathbb{R}_{>0}$ with $B(x,r)\subseteq \Omega$ such that $u\big|_{B(x,r)}$ is Lipschitz. Given a function $u\colon \Omega\to \mathbb{R}$, we define a function $\mathrm{lip}[u]\colon \Omega\to [0,+\infty]$, known as the local Lipschitz constant of $u$, by
\begin{gather}
\mathrm{lip}[u](x)\coloneqq \begin{dcases} \lims\limits_{x'\rightarrow x}\frac{\mathsf{d}_u(x,x')}{\mathsf{d}(x,x')},\quad \text{$x$ is a limit point of $\Omega$},\notag\\
0, \quad\text{otherwise}.
\end{dcases}
\end{gather}
As it can be shown by arguments similar to those provided in \cite[Lemma 6.2.5]{HKST15}, it holds, for any $u\in \mathrm{Lip}_{Loc}(\Omega)$, that $\mathrm{lip}[u]\in \mathfrak{B}(\Omega, \mathbb{R}_{\geq 0})$.

Given a Borel set $E\subseteq \mathrm{X}$, for any $u\in \Leb^1(E)$ we put
\begin{gather}
\langle u\rangle_E\coloneqq \fint\limits_{E} u(x')\D\mathfrak{m}(x')\coloneqq
 \begin{dcases}\tfrac{1}{\mathfrak{m}(E)}\int\limits_{E} u(x')\D\mathfrak{m}(x'),\quad\text{$\mathfrak{m}(E)\in \mathbb{R}_{>0}$},\notag\\
 0,\quad\text{otherwise}\notag.
\end{dcases}
\end{gather}
An analogous notation will also be used several times when integrating functions over subsets of Euclidean spaces with respect to the corresponding Lebesgue measures.

Let $\supp(\mathfrak{m})$ denote the support of $\mathfrak{m}$. Given an open set $\Omega\subseteq \mathrm{X}$, a point $x\in \big(\Omega\cap \supp(\mathfrak{m})\big)$ is called a Lebesgue point of $u\in \Leb^1_{Loc}(\Omega)$ if
\begin{equation*}
\lim\limits_{r\searrow 0}\fint\limits_{B(x,r)}\mathsf{d}_u(x,x') \D \mathfrak{m}(x')=0.
\end{equation*}
Recall also that a point $x\in \supp(\mathfrak{m})$ is called a density point of a Borel set $E\subseteq \mathrm{X}$ if
\begin{equation*}
\lim\limits_{r\searrow 0} \frac{\mathfrak{m}\big(E\cap B(x,r)\big)}{\mathfrak{m}\big(B(x,r)\big)}=1.
\end{equation*}
\subsection{Elements of Sobolev calculus}\label{ss: SobCal}
Now we want to discuss briefly some of those notions that are related to Sobolev calculus on metric measure spaces.

\begin{rema} In his seminal paper \cite{C99}, J. Cheeger developed a theory of Sobolev functions on metric measure spaces. As a part of this, it was proposed how to generalize energies that appear in \eqref{eq: Seminorm} to the singular setting. The corresponding quantities are now occasionally called Cheeger energies, and it is these that we adopt as energies throughout the paper. For a review of this notion in its modern form, we refer the reader to \cite{AD14, AGS13, AGS14}. Several different approaches to the definition of Cheeger energies are provided there. At the same time, as the main results of these papers demonstrate, they turn out to be fully equivalent, at least on complete metric measure spaces.
\end{rema}

We start by defining Cheeger energies for functions. Among possible ways for representing them, in order to describe the cases $p=1$ and $p\in \mathbb{R}_{>1}$ in a unified way, we choose the definition that is based on a relaxation procedure with locally Lipschitz functions involved. Such an approach was proposed for the first time in \cite{M03} by M. Miranda.
\begin{defi}\label{def: ChEnerg}
Let $p\in \mathbb{R}_{\geq 1}$, let $\Omega\subseteq \mathrm{X}$ be an open set, let $u\in \Leb^p_{Loc}(\Omega)$. Given an open set $\Omega'\subseteq \Omega$, we define \textbf{the Cheeger $p$-energy of $u$ on $\Omega'$} as
\begin{equation}
\mathrm{E}_p[u](\Omega')\coloneqq \inf\limits_{(u_n)_{n\in \mathbb{N}}}\limi\limits_{n\rightarrow +\infty} \int\limits_{\Omega'}\big(\mathrm{lip}[u_n](x)\big)^p\D\mathfrak{m}(x),
\label{eq: ChEnerg}
\end{equation}
where the infimum is taken over all sequences $(u_n)_{n\in \mathbb{N}}\subseteq\mathrm{Lip}_{Loc}(\Omega')$ converging to $u$ in $\Leb^p_{Loc}(\Omega')$ as $n\rightarrow +\infty$.
\end{defi}
\noindent In a standard way, we extend the set function $\mathrm{E}_p[u](\bdot)$, which is so far defined only on the family of open subsets of $\Omega $, to an arbitrary set $S\subseteq \Omega$ by letting
\begin{equation*}
\mathrm{E}_p[u](S)\coloneqq \inf\limits_{\Omega'}\mathrm{E}_p[u](\Omega'),
\end{equation*}
where the infimum is taken over all open sets $\Omega'\subseteq \Omega$ with $S\subseteq \Omega'$.

As it turns out, the set function obtained with the above method is a well-posed measure.
\begin{theo}\label{theo: ChAdd}
Let $p\in \mathbb{R}_{\geq 1}$, let $\Omega\subseteq \mathrm{X}$ be an open set, let $u\in\Leb^p_{Loc}(\Omega)$. Then $\mathrm{E}_p[u](\bdot)$ is a Borel-regular outer measure on $\Omega$.
\begin{proof}
For the case $p=1$, the proof can be found in \cite[Theorem 3.4]{M03}. The same arguments, with corresponding minor modifications, are actually applicable to the case $p>1$ as well.
\end{proof}
\end{theo}

We are ready to provide the definition of Cheeger energies for metric-valued maps. It is based on considering post-compositions of a given map with $1$-Lipschitz functions, which reduces everything to the real-valued case. Being proposed by L. Ambrosio in \cite{A90} and Y. G. Reshetnyak in \cite{R97}, this approach appears later in \cite{BNP23,GPS20,GT21}.
\begin{defi}\label{def: MetCh}
Let $p\in \mathbb{R}_{\geq 1}$, let $\Omega\subseteq\mathrm{X}$ be an open set, let $(\mathrm{Y}, \mathsf{d}_{\mathrm{Y}})$ be a metric space, let $f\in \mathfrak{B}(\Omega,\mathrm{Y})$. We define {\textbf{the Cheeger $p$-energy of $f$}}, which we denote by $\mathrm{E}_p[f]$, as the lowest element, in the sense of the lattice of Borel-regular outer measures on $\Omega$, among those $\nu$ that fulfill the following:
\begin{equation*}
\mathrm{E}_p[\phi\circ f](S)\leq \nu(S) \quad \text{for all $S\subseteq \Omega$ and $\phi\in \mathrm{BLip}_1(\mathrm{Y})$}.
\end{equation*}
\end{defi}
\noindent The above definition is well posed, since the considered lattice is known to be complete, and hence the stated supremum always exists. Some related information on the notions involved can be found in \cite{A90} and \cite{HKST15}.

\begin{rema} In contrast to other sources, our definition includes post-compositions with only bounded, rather than arbitrary, $1$-Lipschitz functions. This allows one to define Cheeger energies for all maps in a single manner. A somewhat similar modification appears in \cite[Section 6]{AGS13}.
\end{rema}

The, rather technical, proposition below will be used later.
\begin{prop}\label{prop: SupM}
Let $p\in \mathbb{R}_{\geq 1}$, let $\Omega\subseteq\mathrm{X}$ be an open set, let $(\mathrm{Y}, \mathsf{d}_{\mathrm{Y}})$ be a metric space, let $f\in \mathfrak{B}(\Omega,\mathrm{Y})$. Then the following assertions hold:
\begin{itemize}
\item[$\rm i)$] we have
\begin{equation*}
\mathrm{E}_p[f](\Omega)=\sup\limits_{\Omega'}\mathrm{E}_p[f](\Omega'),
\label{eq: SupM0}
\end{equation*}
where the supremum is taken over all open sets $\Omega'\Subset \Omega$;
\item[$\rm ii)$] we have
\begin{equation*}
\mathrm{E}_p[f](\Omega)=\sup_{(\Omega_j,\phi_j)_{j\in \mathbb{N}}} \sum\limits_{j=1}^{+\infty} \mathrm{E}_p[\phi_j\circ f](\Omega_j),
\label{eq: SupM}
\end{equation*}
where the supremum is taken over all families $(\Omega_j,\phi_j)_{j\in \mathbb{N}}$, where, for any $j\in \mathbb{N}$, $\Omega_j$ is an open subset of $\Omega$, while $\phi_j$ is a function belonging to $\mathrm{BLip}_1(\mathrm{Y})$, such that $(\Omega_j)_{j\in \mathbb{N}}$ is pairwise disjoint.
\end{itemize}
\begin{proof}
Assertion i) follows directly from the additivity property of $\mathrm{E}_p[f]$. For the proof of assertion ii), one should apply some pretty standard measure-theoretic methods.
\end{proof}
\end{prop}

To simplify future calculations, it is useful to introduce the following notation. Let $p\in \mathbb{R}_{\geq 1}$, let $\Omega\subseteq \mathrm{X}$ be an open set, let $(\mathrm{Y}, \mathsf{d}_{\mathrm{Y}})$ be a metric space, let $f\in \mathfrak{B}(\Omega,\mathrm{Y})$. Given $x\in \Omega$ and $r\in \mathbb{R}_{>0}$ with $B(x,r)\subseteq \Omega$, we put
\begin{gather}
\mathcal{A}_p[f](x,r)\coloneqq \begin{dcases}\frac{\mathrm{E}_p[f]\big(B(x, r)\big)}{\mathfrak{m}\big(B(x,r)\big)},\quad \mathfrak{m}\big(B(x,r)\big)\in \mathbb{R}_{>0},\label{eq: AOp}\\
0, \quad \text{otherwise},
\end{dcases}\\
\mathcal{R}_p[f](x,r)\coloneqq \tfrac{1}{3}\sum\limits_{k=0}^{+\infty}\big(\tfrac{2}{3}\big)^k\mathcal{A}_p[f]\big(x,\tfrac{r}{2^k}\big)\label{eq: ROp},\\
\mathcal{M}_p[f](x,r)\coloneqq \sup\limits_{r'\in (0,r]}\mathcal{A}_p[f](x,r')\label{eq: MOp}.
\end{gather}
The just defined quantities can be viewed as, the value at $x$ of, the averaging operator, the Riesz potential, or at least its modified version, and the maximal function, respectively, at scale $r$ applied to $\mathrm{E}_p[f]$. Such regularization operators will be helpful for analyzing the behavior of maps. The reason we defined the Riesz potential in somewhat unusual way will be clarified later. We note that the functions defined via \eqref{eq: AOp}-\eqref{eq: MOp} become Borel whenever $\mathrm{E}_p[f](\Omega)<+\infty$. This can be deduced via arguments provided in \cite[Section 3]{HKST15}.

\begin{rema}
For details on everything said below, we send the reader to \cite{AGS13,AGS14,AD14}. Let $p\in \mathbb{R}_{\geq 1}$, let $\Omega\subseteq \mathrm{X}$ be an open set, let $u\in \Leb^p(\Omega)$. In the case $p>1$, the condition $\mathrm{E}_p[u](\Omega)<+\infty$ can be taken, among other equivalent ones, as a defining property for Sobolev functions. At the same time, in the case $p=1$, this condition defines the family of functions of bounded variation, which only contains that of Sobolev functions as a, typically proper, subfamily. Thereby, other approaches to defining this Sobolev space need to be considered. However, several natural definitions of this space turn out to be non-equivalent. Among those that seem to be satisfactory, probably the most general one is based on the requirement that $\mathrm{E}_1[u](\Omega)<+\infty$ and that $\mathrm{E}_1[u]$ is absolutely continuous with respect to $\mathfrak{m}\resmes_{\Omega}$. But once we take any $p>1$ again, the latter condition automatically implies the former one, so such a definition fits this case as well. And it appears reasonable to expect that all the same should be valid in the metric-valued case.
\end{rema}

Keeping in mind the above remark, we provide now the notion of metric-valued Sobolev maps.
\begin{defi}\label{def: RelGr}
Let $p\in \mathbb{R}_{\geq 1}$, let $\Omega\subseteq \mathrm{X}$ be an open set, let $(\mathrm{Y}, \mathsf{d}_{\mathrm{Y}})$ be a metric space. We define {\rm\textbf{the Sobolev class $\mathrm{S}^{1,p}(\Omega, \mathrm{Y})$}} as the set of all maps $f\in \mathfrak{B}(\Omega, \mathrm{Y})$ with $\mathrm{E}_p[f](\Omega)<+\infty$ such that $\mathrm{E}_p[f]$ is absolutely continuous with respect to $\mathfrak{m}\resmes_\Omega$. The corresponding Radon-Nikodym derivative, powered to the exponent $\frac{1}{p}$, is called {\rm\textbf{the minimal $p$-weak upper gradient of $f$}} and is denoted by $|\nabla f|_p$. We define also {\rm\textbf{the Sobolev space}} as $\mathrm{W}^{1,p}(\Omega, \mathrm{Y})\coloneqq \Leb^p(\Omega, \mathrm{Y})\cap\mathrm{S}^{1,p}(\Omega, \mathrm{Y})$.
\end{defi}

It will be necessary for us to impose a certain extension property on open subsets of $\mathrm{X}$. Not to deal with the highly complicated extension problem in the metric-valued case, we consider extensions with the target space being possibly a larger metric space containing an isometric copy of the original one. And while it obviously suffices to take only Banach spaces as such targets, we will not specify this in our definition. Thus, let $S\subseteq \mathrm{X}$, let $(\mathrm{Y}, \mathsf{d}_{\mathrm{Y}})$ be a metric space, let $f\colon S\to\mathrm{Y}$ be a map, let $(\mathrm{Z},\mathsf{d}_{\mathrm{Z}})$ be another metric space. We say that a map $F\colon \mathrm{X}\to \mathrm{Z}$ is a $\mathrm{Z}$-extension of $f$ if there exists an isometric embedding $\iota\colon \mathrm{Y}\to \mathrm{Z}$ such that $F(x)=\iota\big(f(x)\big)$ for $\mathfrak{m}$-a.e. $x\in S$. Using the introduced notion, we provide the following definition, which extends an analogous one from \cite{LPZ22} to the metric-valued setting.
\begin{defi}\label{def: ExrProp}
Let $p\in \mathbb{R}_{\geq 1}$. We say that an open set $\Omega\subseteq \mathrm{X}$ has {\rm\textbf{the strong $p$-extension property}} if, given an arbitrary metric space $(\mathrm{Y},\mathsf{d}_{\mathrm{Y}})$, any map $f\in \Leb^p(\Omega, \mathrm{Y}\big)$ satisfying the condition $\mathrm{E}_p[f](\Omega)<+\infty$ has, for some metric space $(\mathrm{Z}, \mathsf{d}_{\mathrm{Z}})$, a $\mathrm{Z}$-extension $F$ with $\mathrm{E}_p[F](\mathrm{X})<+\infty$ belonging to $\Leb^p(\mathrm{X}, \mathrm{Z})$ such that
\begin{equation}
\lim\limits_{r\searrow 0} \mathrm{E}_p[F]\big(B(\Omega,r)\big)=\mathrm{E}_p[F](\Omega).
\label{eq: ZerEn}
\end{equation}
\end{defi}

\begin{rema}
For subsequent reasoning, it actually suffices to use a weaker version of the introduced property. Namely, one can require extensions to exist only for maps with the fixed source space. But even with such a relaxation, it is far from obvious whether this property is implied by the classical, i.e. real-valued, strong extension property. Still, we believe that it is supposed to be so for practically important cases, so we do not discuss this issue in more detail.
\end{rema}

Now, following \cite{GT21}, we want to introduce a Lipschitz-like condition, which makes sense to consider for Sobolev maps.
\begin{defi}\label{def: LuzLip}
Let $E\subseteq \mathrm{X}$ be a Borel set, let $(\mathrm{Y}, \mathsf{d}_{\mathrm{Y}})$ be a metric space. We say that a map $f\colon E\to \mathrm{Y}$ has {\rm\textbf{the Luzin-Lipschitz property}} if one can find a family $(E_m)_{m\in \mathbb{N}}$ of Borel subsets of $E$ with $\mathfrak{m}\bigg(E\Big\backslash \bigcup\limits_{m\in \mathbb{N}} E_m\bigg)=0$ such that $f\big|_{E_m}$ is Lipschitz for any $m\in \mathbb{N}$.
\end{defi}

Let us introduce some useful notation, partially adopted from \cite{GT21}. Let $p\in \mathbb{R}_{\geq 1}$, let $\Omega\subseteq \mathrm{X}$ be an open set, let $(\mathrm{Y},\mathsf{d}_{\mathrm{Y}})$ be an metric space, let $f\in \mathfrak{B}(\Omega, \mathrm{Y})$, let $E\subseteq \Omega$ be a Borel set. For all $x\in \Omega$ and $r\in \mathbb{R}_{>0}$, we put
\begin{equation}
\mathsf{ks}_p[f,E](x,r)\coloneqq \begin{dcases}\frac{1}{\mathfrak{m}\big(B(x,r)\big)}\int\limits_{B(x,r)\cap E}\frac{\big(\mathsf{d}_f(x,x')\big)^p}{r^p}\D \mathfrak{m}(x'), \quad \mathfrak{m}\big(B(x,r)\big)\in \mathbb{R}_{>0},\\
0,\quad\text{otherwise}.
\end{dcases}
\label{eq: REnD}
\end{equation}
Such quantities are occasionally called Korevaar-Schoen approximate $p$-energy densities of $f$, after the paper \cite{KS93} on metric-valued Sobolev theory by N. J. Korevaar and R. M. Shoen. They are of interest to us as closely related to quantities as in \eqref{eq: IntSub} but being moreover independent of mollifiers. In terms of them, we define the notion of energy densities, which will be convenient for using thereafter.
\begin{defi}\label{def: EnD}
Let $p\in \mathbb{R}_{\geq 1}$, let $\Omega\subseteq \mathrm{X}$ be an open set, let $(\mathrm{Y},\mathsf{d}_{\mathrm{Y}})$ be an metric space, let $f\in \mathfrak{B}(\Omega, \mathrm{Y})$. The, necessarily $\mathfrak{m}$-a.e. uniquely defined, function $\mathsf{e}_p[f]\in \mathfrak{B}(\Omega, \mathbb{R}_{\geq 0})$, as soon as it exists, will be called {\textbf{the $p$-energy density of $f$}} if, for any open set $\Omega'\Subset \Omega$, it holds that $\mathsf{e}_p[f]\in \Leb^1(\Omega')$ and that the functions $\mathsf{ks}_p[f,\Omega](\bdot,r)$, $r\in \mathbb{R}_{>0}$, converge to $\mathsf{e}_p[f]$ in $\Leb^1(\Omega')$ and $\mathfrak{m}$-a.e. on $\Omega'$ as $r\searrow 0$.
\end{defi}
\subsection{Doubling and Poincar{\'e} conditions}\label{ss: DoubPon}
In this section, we define doubling and Poincar{\'e} conditions, as well as discuss their consequences, that we shall exploit further.

We start by introducing the following doubling property, which can be viewed as a modified version of the semiuniformly local one proposed in \cite{BB18}.
\begin{defi}\label{def: DoubCon}
We say that $\mathrm{X}$ is {\rm\textbf{internally doubling}} if there exists a constant $C_D\in \mathbb{R}_{\geq 1}$ and, for each $S\Subset \mathrm{X}$, there is $R\coloneqq R_D(S)\in \mathbb{R}_{>0}$ such that
\begin{equation*}\label{eq: DoubCon}
0<\mathfrak{m}\big(B(x,2r)\big)\leq C_D\mathfrak{m}\big(B(x,r)\big)<+\infty\quad \text{for all $x\in B(S,R)$ and $r\in (0,R]$}.
\end{equation*}
\end{defi}
\noindent Everywhere below, we will refer to the above relation as the doubling inequality.
\begin{rema}
Loosely speaking, the defined doubling condition means that every well-contained subset of the space has its own doubling radius, whereas the doubling constant should be unique for the whole space, so it would be probably more consistent to call it internally uniformly local one. And in light of that we want to deal mostly with well-contained subsets of our space, such a condition seems to be pretty natural. As a clarification for the reader, we note that a typical case in which the introduced condition may be satisfied, whereas the globally doubling one may not be, is a weighted Euclidean domain. Consider, for instance, the domain $(0,1)$ with the weight $w(x)\coloneqq \tfrac{1}{x}$ and the domain $\mathbb{R}$ with the weight $w(x)\coloneqq \exp(x)$. For some additional information on various local doubling conditions and where they usually appear, we refer the reader to \cite{BB18}.
\end{rema}

In the following proposition, we list the geometric properties of the underlying space that are implied by the just introduced doubling condition.
\begin{prop}\label{prop: LocCom}
Suppose $\mathrm{X}$ is internally doubling. Then we have $\supp(\mathfrak{m})=\mathrm{X}$. Moreover, for all $S\Subset \mathrm{X}$ and $x\in S$, the ball $B(x,R)$, where $R\coloneqq \tfrac{1}{2}\min\Big\{\dist\big(S, \overline{\mathrm{X}}\backslash \mathrm{X}\big), R_D(S)\Big\}$, is compact.
\begin{proof}
The first part follows directly from \cref{def: DoubCon}. A possible proof for the second part can be adopted from \cite[Section 4.1]{HKST15}.
\end{proof}
\end{prop}

The proposition below concerns the standard fact that, under doubling conditions, the Lebesgue differentiation theorem is valid.
\begin{prop}\label{prop: LebDif}
Suppose $\mathrm{X}$ is internally doubling. Let $\Omega\subseteq \mathrm{X}$ be an open set, let $u\in \Leb^1_{Loc}(\Omega)$. Then $\mathfrak{m}$-a.e. point $x\in \Omega$ is a Lebegsue point of $u$.
\begin{proof}
For the proof, we send the reader to \cite[Section 3.4]{HKST15}.
\end{proof}
\end{prop}

In addition to the proposition above, our doubling condition yields the following properties for the transformations defined via \eqref{eq: AOp}-\eqref{eq: MOp}.
\begin{prop}\label{prop: RegProp}
Suppose $\mathrm{X}$ is internally doubling. Let $p\in \mathbb{R}_{\geq 1}$, let $\Omega\subseteq \mathrm{X}$ be an open set, let $(\mathrm{Y},\mathsf{d}_{\mathrm{Y}})$ be a metric space, let $f\in \mathfrak{B}(\Omega,\mathrm{Y})$. Fix a Borel set $E\Subset \Omega$ and put $R\coloneqq \tfrac{1}{16}\min\Big\{\dist\big(E,\overline{\mathrm{X}}\backslash \Omega\big), R_D(E)\Big\}$. Then the following assertions are satisfied:
\begin{itemize}
\item[$\rm i)$] if $\mathrm{E}_p[f](\Omega)<+\infty$, then for any $r\in (0,R]$ we have
\begin{equation*}
\max\biggg\{\int\limits_{E}\mathcal{A}_p[f](x,r)\D \mathfrak{m}(x),\int\limits_{E}\mathcal{R}_p[f](x,r)\D \mathfrak{m}(x)\biggg\}\leq C_D\mathrm{E}_p[f]\big(B(E,r)\big);
\end{equation*}
\item[$\rm ii)$] if $f\in \mathrm{S}^{1,p}(\Omega,\mathrm{Y})$, then $\big(\mathcal{A}_p[f](\bdot,r)\big)_{r\in (0,R)}$ and $\big(\mathcal{R}_p[f](\bdot,r)\big)_{r\in (0,R)}$ converge to $|\nabla f|_p^p$ in $\Leb^1(E)$ as $r\searrow 0$;
\item[$\rm iii)$] if $\mathrm{E}_p[f](\Omega)<+\infty$, it follows that $\mathcal{M}_p[f](\bdot,r)$ is $\mathfrak{m}$-a.e. finite on $E$ for any $r\in (0,R]$.
\end{itemize}
\begin{proof}
The parts of assertions i) and ii) about $\mathcal{A}_p[f]$ can be established in a way similar to that given in \cite[Lemma 3.11]{GT21}. They in turn imply the corresponding parts about $\mathcal{R}_p[f]$, by the very definitions of both. Assertion iii) is a trivial consequence of the well-known weak-type estimates applied to $\mathrm{E}_p[f]$, see \cite[Theorem 3.5.6]{HKST15} for a possible reasoning. We note that in the sources we mentioned, stronger doubling conditions, compared to the one we use, are required for the space. Still, since we consider sufficiently small radii only, we insist that the proofs from the sources can be applied to our case as well.
\end{proof}
\end{prop}

One more well-known fact is that one can construct a Lipschitz partition of unity with certain useful properties on spaces satisfying an appropriate doubling condition. We shall use this fact in the following formulation.
\begin{prop}\label{prop: PartUnit}
Suppose $\mathrm{X}$ is internally doubling. Then there exists a constant $c_{PU}\in \mathbb{R}_{>0}$ depending only on $C_D$ such that the following holds. For each $S\Subset \mathrm{X}$ there is $r_{PU}\in \mathbb{R}_{>0}$ such that, for any $r\in(0,r_{PU})$, one can find a data set that consists of points $x_i\in S$, $i\in \mathbb{N}$, and of functions $\varphi_i\colon \mathrm{X}\to  [0,1]$, $i\in \mathbb{N}$, with the following properties:
\begin{itemize}
\item[$\rm i)$] for every $x\in \mathrm{X}$ we have $\sum\limits_{i\in \mathbb{N}} \chi_{B(x_i,2r)}(x)\leq c_{PU}$;
\item[$\rm ii)$] for every $i\in \mathbb{N}$ it follows that $\varphi_i$ is $\tfrac{c_{PU}}{r}$-Lipschitz and that $\varphi_i(x)=0$ for any $x\in \mathrm{X}\backslash B(x_i,r)$;
\item[$\rm iii)$] for every $x\in S$ we have $\sum\limits_{i\in \mathbb{N}} \varphi_i(x)=1$.
\end{itemize}
\begin{proof}
For the proof, we refer to \cite[Lemma 2.4]{GT21}, where a similar statement is proved for the whole space, not for its subset, and under the assumption that the source space satisfies a stronger doubling condition than the one we use. However, since $S$ is a well-contained subset of $\mathrm{X}$ and since we claim the existence of a partition of unity that is subordinate to a cover consisting of balls of sufficiently small radii only, the proof of the mentioned lemma can be implemented for our case with some trivial minor modifications.
\end{proof}
\end{prop}

Now we introduce a Poincar{\'e} condition that will be considered in what follows.
\begin{defi}\label{def: PoinIn}
Let $p\in \mathbb{R}_{\geq 1}$. We say that $\mathrm{X}$ supports {\rm\textbf{an internal $p$-Poincar{\'e} inequality}} if there exist constants $C_P\in \mathbb{R}_{>0}$ and $\lambda\coloneqq \lambda_P\in \mathbb{R}_{\geq 1}$ and, for each $S\Subset \mathrm{X}$, there is $R\coloneqq R_P(S)\in \mathbb{R}_{>0}$ such that $0<\mathfrak{m}\big(B(x,r)\big)<+\infty$ for all $x\in B(S,R)$ and $r\in (0,\lambda R]$, as well as
\begin{equation*}
\fint\limits_{B(x,r)}\big|u(x')-\langle u\rangle_{B(x,r)}\big|\D\mathfrak{m}(x')\leq C_P r \Biggg(\fint\limits_{B(x,\lambda r)}\big(\mathrm{lip}[u](x')\big)^p\D \mathfrak{m}(x')\Biggg)^\frac{1}{p}
\end{equation*}
for any Lipschitz function $u\colon \mathrm{X}\to \mathbb{R}$ and for all $x\in B(\Omega, R)$ and $r\in (0,R]$.
\end{defi}
\noindent Everything that is said within the remark under \cref{def: DoubCon} corresponds to the given Poincar{\'e} condition too.
\begin{rema}
We note that in the literature, for example, in \cite{HKST15}, Poincar{\'e} conditions are typically formulated with the use of so-called upper gradients. And while the condition from the definition above is a priori less restrictive, results provided in \cite{AD14, AGS13} demonstrate that both variants are perfectly equivalent.
\end{rema}

As usual, a combination of \cref{def: ChEnerg} and \cref{def: PoinIn} yields an upgraded Poincar{\'e} inequality. We state the corresponding fact, for which we also recall the notation given in \eqref{eq: AOp}, in the proposition below.
\begin{prop}\label{prop: PoinIn}
Let $p\in \mathbb{R}_{\geq 1}$. Suppose $\mathrm{X}$ supports an internal $p$-Poincar{\'e} inequality. Let $\Omega\subseteq \mathrm{X}$ be an open set, let $u\in \Leb^p_{Loc}(\Omega)$. Fix $S\Subset \Omega$, put $\lambda \coloneqq \lambda_P$ and $R\coloneqq R_P(S)$. Then we have
\begin{equation*}
\fint\limits_{B(x,r)}\big|u(x')-\langle u\rangle_{B(x,r)}\big|\D\mathfrak{m}(x')\leq C_P r\Big(\mathcal{A}_p[u](x,\lambda r)\Big)^{\frac{1}{p}} \quad \text{for all $x\in B(S,R)$ and $r\in (0,R]$}
\label{eq: PoinIn}
\end{equation*}
whenever the ball $B(x,\lambda r)$ is a compact subset of $\Omega$.
\begin{proof}
The proof is standard, so we omit it. Some details can be found in \cite[Lemma 1.13]{BPR20}.
\end{proof}
\end{prop}
\noindent We will refer to the inequality above as the Poincar{\'e} inequality.
\subsection{Requirements for mollifiers}\label{ss: Mol}
This section is aimed at indicating assumptions on mollifiers that we want to deal with.

Those conditions for mollifiers that we consider throughout the article are listed in the following definition. Their choice is determined by a method that we shall exploit in further proofs for representing quantities as in \eqref{eq: IntSub} in terms of quantities as in \eqref{eq: REnD}. They carry a somewhat similar meaning to those used in \cite{LPZ22} and \cite{HP21}, even though they are more cumbersome.
\begin{defi}\label{def: Mol}
Let $p\in \mathbb{R}_{\geq 1}$, let $(\rho_{\delta})_{\delta\in (0,1)}\subseteq \mathfrak{B}\big(\mathrm{X}\times \mathrm{X}, \mathbb{R}_{\geq 0}\big)$. We say that the family $(\rho_{\delta})_{\delta\in (0,1)}$ is \textbf{$p$-admissible} if the following conditions are satisfied:
\begin{itemize}
\item[$\rm i)$] it holds that
\begin{gather}
\lims\limits_{r\searrow 0}\lims\limits_{\delta \searrow 0} \sup\limits_{x\in \mathrm{X}} \int\limits_{\mathrm{X}\backslash B(x,r)} \big(\rho_{\delta}(x,x')+\rho_{\delta}(x',x)\big) \D \mathfrak{m}(x')=0;\label{eq: Cond1}
\end{gather}
\item[$\rm ii)$] one has
\begin{gather}
\lims\limits_{r\searrow 0}\lims\limits_{\delta\searrow 0} \sigma(\delta,x,r)=0\quad \text{for any $x\in \mathrm{X}$}\label{eq: Cond2},
\end{gather}
where, for all $x\in \mathrm{X}$, $\delta\in (0,1)$, and $r\in \mathbb{R}_{>0}$, we put
\begin{equation}
\sigma(\delta,x,r)\coloneqq \begin{dcases}
\sup\Big\{\rho_{\delta}(x,x')\bigm| x'\in \big(\mathrm{X}\backslash B(x,r)\big)\Big\},\quad \big(\mathrm{X}\backslash B(x,r)\big)\neq \emptyset,\\
0, \quad \text{otherwise};
\end{dcases}\label{eq: MNotS}
\end{equation}
\item[$\rm iii)$] for all $\delta\in (0,1)$ and $x,x_1,x_2\in \mathrm{X}$ with $\mathsf{d}(x,x_1)\leq \mathsf{d}(x,x_2)$ one has $\rho_{\delta}(x,x_1)\geq \rho_{\delta}(x,x_2)$;
\item[$\rm iv)$] one has
\begin{equation}
\lims\limits_{r\searrow 0}\Big(\sigma(\delta,x,r)\mathfrak{m}\big(B(x,r)\big)r^p\Big)=0\quad\text{for all $x\in \mathrm{X}$ and $\delta\in (0,1)$};
\label{eq: Cond4}
\end{equation}
\item[$\rm v)$] there is a constant $C_M\in \mathbb{R}_{\geq 1}$  that fulfills, for any nonempty set $S\Subset \mathrm{X}$, the inequalities
\begin{gather}
\frac{1}{C_M}\leq \limi\limits_{r\searrow 0}\limi\limits_{\delta\searrow 0}\sum\limits_{k=0}^{+\infty} \inf\limits_{x\in S}\Pi_{L,k}(\delta, x,r,2)\leq \lims\limits_{r\searrow 0}\lims\limits_{\delta\searrow 0}\sum\limits_{k=0}^{+\infty} \sup\limits_{x\in S}\Pi_{U,k}(\delta, x,r,2)\leq C_M,
\label{eq: Cond5}
\end{gather}
where we put
\begin{gather}
\Pi_{L,k}(\delta,x,r,h)\coloneqq \begin{dcases}\Big|\sigma\big(\delta,x,\tfrac{r}{h^{k}}\big)-\sigma\big(\delta,x,\tfrac{r}{h^{k-1}}\big)\Big|\mathfrak{m}\Big(B\big(x,\tfrac{r}{h^{k+1}}\big)\Big)\big(\tfrac{r}{h^{k+1}}\big)^p,\quad k\neq 0,\\
\sigma(\delta,x,r)\mathfrak{m}\Big(B\big(x,\tfrac{r}{h}\big)\Big)\big(\tfrac{r}{h}\big)^p,\quad \text{otherwise},
\end{dcases} \label{eq: MNotPL}\\
\Pi_{U,k}(\delta,x,r,h)\coloneqq \begin{dcases}\Big|\sigma\big(\delta,x,\tfrac{r}{h^{k}}\big)-\sigma\big(\delta,x,\tfrac{r}{h^{k-1}}\big)\Big|\mathfrak{m}\Big(B\big(x,\tfrac{r}{h^{k}}\big)\Big)\big(\tfrac{r}{h^{k}}\big)^p,\quad k\neq 0,\\
\sigma(\delta,x,r)\mathfrak{m}\big(B(x,r)\big)r^p,\quad \text{otherwise},
\end{dcases}\label{eq: MNotPU}
\end{gather}
for any $k\in \mathbb{N}_0$ and for all $x\in \mathrm{X}$, $\delta\in (0,1)$, $r\in \mathbb{R}_{>0}$, and $h\in \mathbb{R}_{>1}$.
\end{itemize}
We also say that the family $(\rho_{\delta})_{\delta\in (0,1)}$ is \textbf{strongly $p$-admissible} if, in addition to conditions i)-v) above, it satisfies the following ones:
\begin{itemize}
\item[$\rm vi)$] one has
\begin{equation}
\lims\limits_{r\searrow 0}\lims\limits_{\delta\searrow 0}\sum\limits_{k=0}^{+\infty} \sup\limits_{x\in S}\Pi_{U,k}(\delta, x,r,h)<+\infty \quad \text{for all $S\Subset \mathrm{X}$ and $h\in \mathbb{R}_{>1}$};
\label{eq: Cond6}
\end{equation}
\item[$\rm vii)$] there exists a function $\Theta_M\in \mathfrak{B}(\mathrm{X}, \mathbb{R}_{\geq 0})$ such that
\begin{gather}
\lims\limits_{h\searrow 1}\esup\limits_{x\in \mathrm{X}}\lims\limits_{r\searrow 0}\lims\limits_{\delta\searrow 0}\biggg(\Bigg|\Theta_M(x)-\sum\limits_{k=0}^{+\infty}\Pi_{L,k}(\delta,x,r,h)\Bigg|+\Bigg|\Theta_M(x)-\sum\limits_{k=0}^{+\infty}\Pi_{U,k}(\delta,x,r,h)\Bigg|\biggg)=0.
\label{eq: Cond7}
\end{gather}
\end{itemize}
\end{defi}
\noindent The conditions given above will be called conditions i)-vii) of mollifiers later. Note that the quantities defined via \eqref{eq: MNotS} are always finite as our space is connected. In view of this, the quantities given by \eqref{eq: MNotPL} and \eqref{eq: MNotPU}, in turn, are well-posed and finite at sufficiently small scales.

\begin{rema}
Let us give some comments on \cref{def: Mol}. 

First, we briefly clarify the essence of the listed conditions. Conditions i) and ii) mean that mollifiers become increasingly localized as $\delta\searrow 0$. It will allow all issues related to mollifiers to be considered locally. Condition iii) states that mollifiers are monotonically non-increasing with respect to the distance between points provided that the first one is fixed, which is naturally to expect from objects of this type. Condition iv) imposes a restriction on the asymptotic behavior of mollifiers as $r\searrow 0$. It is a necessary requirement for our subsequent arguments to be applicable. Within condition v), final estimates that are directly used in the proof of \cref{theo: Res1} are formulated. Condition vi), which means basically the same as condition v), but for all values of the multiplicative parameter $h$, will enable us to interchange limits and integrals. Finally, condition vii) imposes on mollifiers a strong regularity property that is needed for the proof of \cref{theo: Res2}.

Note that some of the listed conditions can be slightly relaxed even within the same framework. As an example, in light of subsequent proofs, one can declare a given family of mollifiers admissible if it is "equivalent", up to a universal multiplicative constant, to another family of mollifiers that, in turn, satisfies conditions i)-v). All other possible relaxations can be detected from further reasoning if necessary.

Clearly, an arbitrary metric measure space may have no admissible, much less strongly admissible, families of mollifiers, since their existence implicitly imposes on the underlying space some restrictions. But it will be demonstrated later that this is not the case for \cref{theo: Res1} and \cref{theo: Res2}, i.e. some concrete examples of families of those types that are considered there can actually be granted.
\end{rema}
\subsection{Strongly rectifiable spaces}\label{ss: StrRectSp}
Our goal for this section is to formulate one highly important result about strongly rectifiable spaces.

Informally speaking, strongly rectifiable spaces are those metric measure spaces that can, for any given $\varepsilon\in (0,1)$, be partitioned, up to a set of zero measure, into a countably infinite family of Borel sets admitting a $(1+\varepsilon)$-biLipschitz map into a, possibly depending on the set, Euclidean space such that the pushforward measure obtained via this map is approximately, with the multiplicative error being less than $(1+\varepsilon)$, proportional to the corresponding Lebesgue measure. We refer to \cite{GP22} as the original definition and to \cite{GT21} as its slightly modified, but still equivalent, version. We emphasize that, in contrast to \cite{GT21}, we do not fix for a strongly rectifiable space its dimension, i.e. the dimensions of target Euclidean spaces.

It follows directly from the definition of strongly rectifiable spaces that they possess a dimensional decomposition. This fact, in the case when $\mathrm{X}$ is strongly rectifiable, leads to the existence of a function $\mathrm{Dim}\in \mathfrak{B}(\mathrm{X}, \mathbb{N})$ that indicates for a point the dimension of its component. Accordingly, for any given $d\in \mathbb{N}$, let $\mathrm{X}(d)$ denote the corresponding dimensional component, i.e. the preimage of $\{d\}$ under $\mathrm{Dim}$.

\begin{rema}
The main reason we focus exactly on strongly rectifiable spaces is one result about them, obtained in \cite{GT21}. As its first part, it was shown that, when being defined on a strongly rectifiable metric measure space supporting an appropriate doubling condition, metric-valued maps having the Luzin-Lipschitz property in the sense of \cref{def: LuzLip}, are approximately metrically differentiable. This notion extends that of metric differentiability, introduced in \cite{K94} by B. Kirchheim, to the setting singular with respect to the source space. Being a consequence of the first part, the second one states that a large class of metric-valued maps, including Sobolev ones, admit energy densities in a sense analogous to the given in \cref{def: EnD}. This part however is established only in the case $p>1$. This is because its proof is highly based on the characterization of Sobolev functions through Haj{\l}asz gradients. For the corresponding definitions, we refer to \cite{H09} by P. Haj{\l}asz. These connections, in turn, are directly related to the strong-type estimates, which are valid in the case $p>1$ only.
\end{rema}

We want to eliminate the drawback described in the above remark, i.e. to obtain eventually the existence of energy densities for the case $p=1$ as well, so we provide here only that part of the mentioned result that is applicable for any $p\in \mathbb{R}_{\geq 1}$. The formulation is as follows.
\begin{theo}\label{theo: ApprMetrDiff}
Suppose $\mathrm{X}$ is strongly rectifiable. Let $p\in \mathbb{R}_{\geq 1}$. Suppose $\mathrm{X}$ is internally doubling. Let $\Omega\subseteq \mathrm{X}$ be an open set, let $(\mathrm{Y}, \mathsf{d}_{\mathrm{Y}})$ be a metric space, let $f\in \mathfrak{B}(\Omega, \mathrm{Y})$ be a map having the Lusin-Lipschitz property. Then there exists a function $e\in \mathfrak{B}(\Omega,\mathbb{R}_{\geq 0})$ such that, for $\mathfrak{m}$-a.e. $x\in \Omega$, one can find a family $(V_{x,n})_{n\in \mathbb{N}}$ of Borel subsets of $\Omega$ having $x$ as a density point satisfying the following property:
\begin{equation}
\lims\limits_{n\rightarrow +\infty}\lims\limits_{r \searrow 0} \Big|\mathsf{ks}_p[f,V_{x,n}](x,r)-e(x)\Big|=0.
\label{eq: ApprMetrDiff}
\end{equation}
Furthermore, for each $d\in \mathbb{N}$ and for $\mathfrak{m}$-a.e. $x\in \mathrm{X}(d)$ there is a seminorm $\mathsf{md}_x[f]$ on $\mathbb{R}^d$ such that
\begin{equation}
e(x)\coloneqq \fint\limits_{B^d}\big(\mathsf{md}_x[f](x')\big)^p\D \mathcal{L}^d(x').
\label{eq: ApprMetrDiff0}
\end{equation}
\begin{proof}
The assertion follows from the proofs of \cite[Proposition 3.6]{GT21} and \cite[Theorem 3.13]{GT21}, but, for the convenience of the reader, we provide here some explanations. To satisfy completely those assumptions that are used in \cite{GT21}, suppose for a moment that $\mathrm{X}$ satisfies the corresponding, stronger than the our one, doubling condition and that $\mathrm{X}=\mathrm{X}(d)=\Omega$ for some $d\in \mathbb{N}$.

It follows that, for any given sequence $(\varepsilon_n)_{n\in \mathbb{N}}\subseteq (0,1)$ converging to zero as $n\rightarrow +\infty$, one can construct a certain family $\mathfrak{A}\coloneqq\Big(\big(V_{d,n,l}, \psi_{d,n,l}\big)_{l\in \mathbb{N}}\Big)_{n\in \mathbb{N}}$, called an aligned family of atlases, where $V_{d,n,l}$ is a Borel subset of $\mathrm{X}$ and $\psi_{d,n,l}$ is a function from $V_{d,n,l}$ to $\mathbb{R}^d$ for all $n,l\in \mathbb{N}$. One of its properties is that for $\mathfrak{m}$-a.e. $x\in \mathrm{X}$ and for any $n\in \mathbb{N}$ there is $l(x,n)\in \mathbb{N}$ such that $x$ is a density point of $V_{d,n,l(x,n)}$. Next, \cite[Proposition 3.6]{GT21} states that $f$, as having the Luzin-Lipschitz property, is approximately metrically differentiable with respect to $\mathfrak{A}$ at $\mathfrak{m}$-a.e. point $x\in \mathrm{X}$, which means, in particular, that there exists a special seminorm $\mathsf{md}_x[f]$ on $\mathbb{R}^d$, called the approximate metric differential of $f$ with respect to $\mathfrak{A}$ at $x$. From the arguments given within the second part of the proof of \cite[Theorem 3.13]{GT21}, it follows for $\mathfrak{m}$-a.e. $x\in \mathrm{X}$ that
\begin{equation*}
\lims\limits_{n\rightarrow +\infty}\lims\limits_{r \searrow 0} \biggg|\mathsf{ks}_p\big[f,V_{d,n,l(x,n)}\big](x,r)-\fint\limits_{B^d}\big(\mathsf{md}_x[f](x')\big)^p\D \mathcal{L}^d(x')\biggg|=0.
\end{equation*}
The just applied theorem is exactly the place where the authors restrict themselves to the case $p>1$, but we claim however that all the reasoning, only in the part we explained, is still valid in the case $p=1$.

From all that was said, letting $Q$ be the set of all "good" points for which all the discussed properties hold, putting $V_{x,n}\coloneqq V_{d,n,l(x,n)}$ for all $x\in Q$ and $n\in \mathbb{N}$, and defining a function $e\colon \mathrm{X}\to \mathbb{R}_{\geq 0}$ by
\begin{equation*}
e(x)\coloneqq \begin{dcases}
\fint\limits_{B^d}\big(\mathsf{md}_x[f](x')\big)^p\D \mathcal{L}^d(x'),\quad x\in Q,\\
0, \quad \text{otherwise},
\end{dcases}
\label{eq: SizeD}
\end{equation*}
we end up with the desired statement, up to the suppositions we made. But all computations provided in \cite{GT21} are essentially local and can be carried out independently for each dimensional component, so the conclusion of theorem indeed follows.
\end{proof}
\end{theo}

\begin{rema}
We will prove eventually that the above function $e$ is, under certain conditions, the $p$-energy density of $f$ in the sense of \cref{def: EnD}. Thus, taking into account \eqref{eq: ApprMetrDiff0}, we obtain a quite transparent formula for this density, namely in terms of the metric differential of $f$.
\end{rema}

We conclude the section by providing the following theorem that shows that the measure on a strongly rectifiable space behaves regularly at small scales. It will be used afterwards to verify that such spaces possess strongly admissible families of mollifiers in the sense of \cref{def: Mol}.
\begin{theo}\label{theo: HausDens}
Suppose $\mathrm{X}$ is strongly rectifiable. Then for $\mathfrak{m}$-a.e. $x\in \mathrm{X}$ it holds that
\begin{equation}
\lim\limits_{r\searrow 0} \frac{\mathfrak{m}\big(B(x,hr)\big)}{\mathfrak{m}\big(B(x,r)\big)}=h^{\mathrm{Dim}(x)}\quad \text{for every $h\in \mathbb{R}_{>0}$}.
\label{eq: HausDens0}
\end{equation}
\begin{proof}
Although the assertion is an immediate consequence of \cite[Theorem 2.13]{DG18}, let us give some details here. Fix $d\in \mathbb{N}$. Let $\mathcal{H}^d$ be the $d$-dimensional Hausdorff measure on $\mathrm{X}$. Define a function $\theta_d\colon \mathrm{X}(d)\to \mathbb{R}_{\geq 0}$ by
\begin{equation*}
\theta_d(x)\coloneqq \begin{dcases}
\lim\limits_{r\searrow 0} \frac{\mathfrak{m}\big(B(x,r)\big)}{\mathcal{L}^d(B^d) r^d}, \quad \text{the limit exists and is finite},\\
0, \quad \text{otherwise}.
\end{dcases}
\end{equation*}
It is shown in \cite[Theorem 2.13]{DG18} that the function $\theta_d$ is a Borel representative of the Radon-Nikodym derivative of $\mathfrak{m}\resmes_{\mathrm{X}(d)}$ with respect to $\mathcal{H}^d\resmes_{\mathrm{X}(d)}$. It is obvious then that the measure, with respect to $\mathfrak{m}$, of the set of all $x\in \mathrm{X}(d)$ with $\theta_d(x)=0$ is zero. This, together with the definition of $\mathrm{X}(d)$, trivially ensures that \eqref{eq: HausDens0} holds for $\mathfrak{m}$-a.e. $x\in \mathrm{X}(d)$. The conclusion of the assertion follows then from the arbitrariness of $d$.
\end{proof}
\end{theo}
\section{Proofs}\label{ss: Proof}
\subsection{Auxiliary statements}\label{ss: AuxSt} Here we establish several useful technical facts that will allow us to prove \cref{theo: Res1} and \cref{theo: Res2}.

We start with the simple proposition below, which, roughly speaking, will reduce all further reasoning to local issues. For its proof, we recall the conditions from \cref{def: Mol}, as well as \cref{prop: LipSup}.
\begin{prop}\label{lem: BExt}
Let $p\in \mathbb{R}_{\geq 1}$, let $(\rho_{\delta})_{\delta\in (0,1)}$ be a $p$-admissible family. Let $\Omega\subseteq \mathrm{X}$ be an open set, let $(\mathrm{Y}, \mathsf{d}_{\mathrm{Y}})$ be a metric space, let $f\in \Leb^p(\Omega, \mathrm{Y})$. Then the following properties hold:
\begin{gather}
\lims\limits_{r\searrow 0}\lims\limits_{\delta\searrow 0} \int\limits_{\Omega}\int\limits_{\Omega\backslash B(x,r)}\big(\mathsf{d}_f(x,x')\big)^p\rho_{\delta}(x,x')\D \mathfrak{m}(x')\D \mathfrak{m}(x)=0;\label{eq: BExt0}\\
\lims\limits_{r\searrow 0}\lims\limits_{\delta\searrow 0}\int\limits_{\Omega\backslash B(x,r)} \big(\mathsf{d}_f(x,x')\big)^p\rho_{\delta}(x,x')\D \mathfrak{m}(x')=0 \quad\text{for any $x\in \Omega$}.\label{eq: BExt00}
\end{gather}
\begin{proof}
Since $f\in \Leb^p(\Omega,\mathrm{Y})$, it follows that $\mathsf{d}_{\mathrm{Y}}\big(f(\bdot),y\big)\in \Leb^p(\Omega)$ for some $y\in \mathrm{Y}$. Applying Jensen's inequality and Tonneli's theorem, which can be used here due to \cref{prop: LipSup}, one can easily obtain for all $\delta\in (0,1)$ and $r\in \mathbb{R}_{>0}$ that
\begin{gather*}
\int\limits_{\Omega}\int\limits_{\Omega\backslash B(x,r)}\big(\mathsf{d}_f(x,x')\big)^p\rho_{\delta}(x,x')\D \mathfrak{m}(x')\D \mathfrak{m}(x)\leq\notag\\
\leq 2^{p-1}\biggg(\int\limits_{\Omega} \Big(\mathsf{d}_{\mathrm{Y}}\big(f(x),y\big)\Big)^p \D \mathfrak{m}(x)\biggg)\sup\limits_{x\in \Omega} \int\limits_{\Omega\backslash B(x,r)} \big(\rho_{\delta}(x,x
')+\rho_{\delta}(x',x)\big)\D \mathfrak{m}(x').
\end{gather*}
In view of condition i) of mollifiers, after passing to the upper limits as $\delta\searrow 0$ and then as $r\searrow 0$, we come to \eqref{eq: BExt0}. After that, following the same strategy as above and keeping in mind the notation given in \eqref{eq: MNotS}, we deduce, for all $\delta\in (0,1)$, $x\in \Omega$, and $r\in \mathbb{R}_{>0}$, that
\begin{gather*}
\int\limits_{\Omega\backslash B(x,r)} \big(\mathsf{d}_f(x,x')\big)^p\rho_{\delta}(x,x')\D \mathfrak{m}(x')\leq \notag\\
\leq 2^{p-1}\Big(\mathsf{d}_{\mathrm{Y}}\big(f(x),y\big)\Big)^p\int\limits_{\Omega\backslash B(x,r)}\rho_{\delta}(x,x')\D \mathfrak{m}(x')+\notag\\
+ 2^{p-1}\sigma(\delta,x,r)\int\limits_{\Omega\backslash B(x,r)}\Big(\mathsf{d}_{\mathrm{Y}}\big(f(x'),y\big)\Big)^p\D \mathfrak{m}(x').
\end{gather*}
By conditions i) and ii) of mollifiers, after passing to the upper limits as $\delta\searrow 0$ and then as $r\searrow 0$, we get \eqref{eq: BExt00}. This finishes the proof.
\end{proof}
\end{prop}

The next lemma demonstrates our main method, on which the subsequent proofs are based, for establishing estimates between quantities as in \eqref{eq: IntSub} and as in \eqref{eq: REnD} in terms of the objects defined via \eqref{eq: MNotPL} and \eqref{eq: MNotPU}. The proof relies on a careful use of the conditions from \cref{def: Mol}.
\begin{lem}\label{lem: KeyL}
Let $p\in \mathbb{R}_{\geq 1}$, let $(\rho_{\delta})_{\delta\in (0,1)}$ be a $p$-admissible family. Let $\Omega\subseteq \mathrm{X}$ be an open set, let $(\mathrm{Y}, \mathsf{d}_{\mathrm{Y}})$ be a metric space, let $f\in \mathfrak{B}(\Omega, \mathrm{Y})$. Fix $\delta\in (0,1)$, $x\in \mathrm{X}$ and $r\in \mathbb{R}_{>0}$ with $\mathfrak{m}\big(B(x,r)\big)<+\infty$, $h\in \mathbb{R}_{>1}$, and a Borel set $E\subseteq \Omega$. Then we have
\begin{gather}
\sum\limits_{k=0}^{+\infty} \Pi_{L,k}(\delta,x,r,h) \mathsf{ks}_p[f,E]\big(x,\tfrac{r}{h^{k+1}}\big)\leq \notag\\
\leq  \int\limits_{E\cap B(x,r)}\big(\mathsf{d}_f(x,x')\big)^p\rho_{\delta}(x,x') \D \mathfrak{m}(x')\leq \notag\\
\leq \sum\limits_{k=0}^{+\infty}\Pi_{U,k}(\delta,x,r,h)\mathsf{ks}_p[f,E]\big(x,\tfrac{r}{h^k}\big).
\label{eq: KeyL}
\end{gather}
\begin{proof}
Fix $\beta\in \mathbb{R}_{>0}$. Define a function $\Psi_\beta\in \mathfrak{B}\big(B(x,r),\mathbb{R}_{\geq 0}\big)$ by
\begin{equation}
\Psi_{\beta}(x')\coloneqq\min\Big\{\chi_E(x')\big(\mathsf{d}_f(x,x')\big)^p, \beta \big(\mathsf{d}(x,x')\big)^p\Big\}.
\label{eq: KeyL0}
\end{equation} 
For any $k\in \mathbb{N}_0$ put
\begin{equation}
J_k\coloneqq \int\limits_{B\big(x,\frac{r}{h^k}\big)}\chi_E(x')\big(\mathsf{d}_f(x,x')\big)^p\D \mathfrak{m}(x'), \quad J_k(\beta)\coloneqq \int\limits_{B\big(x,\frac{r}{h^k}\big)}\Psi_{\beta}(x')\D \mathfrak{m}(x').
\label{eq: KeyL1}
\end{equation}
Notice that
\begin{equation}
J_k(\beta)\leq \beta \mathfrak{m}\big(B(x,r)\big)r^p<+\infty\quad \text{for any $k\in \mathbb{N}_0$}.
\label{eq: KeyL2}
\end{equation}
Keeping in mind the notation given in \eqref{eq: MNotS}, we put $\sigma_k\coloneqq \sigma\big(\delta,x,\tfrac{r}{h^{k}}\big)$ for any $k\in \mathbb{N}$ and also $\sigma_0\coloneqq 0$. Notice that the sequence $(\sigma_k)_{k\in \mathbb{N}_0}$ is non-decreasing. Consider the partition of the set $\big(B(x,r)\backslash \{x\}\big)$ into the sets $\Big(B\big(x,\tfrac{r}{h^k}\big)\backslash B\big(x,\tfrac{r}{h^{k+1}}\big)\Big)$, $k\in\mathbb{N}_0$. From \eqref{eq: KeyL0}, we deduce that $\Psi_{\beta}(x)=0$. Thus, in view of \eqref{eq: KeyL1} and of condition iii) of mollifiers, we can write
\begin{gather}
\int\limits_{B(x,r)}\Psi_{\beta}(x')\rho_{\delta}(x,x') \D \mathfrak{m}(x')\leq \sum\limits_{k=0}^{+\infty} \sigma_{k+1}\big(J_{k}(\beta)-J_{k+1}(\beta)\big)=\notag\\
=\sum\limits_{k=0}^{+\infty} \bigg(\Big(\sigma_kJ_k(\beta)-\sigma_{k+1}J_{k+1}(\beta)\Big)+J_k(\beta)\big(\sigma_{k+1}-\sigma_k\big)\bigg).
\label{eq: KeyL4}
\end{gather}
By combining condition iv) of mollifiers with \eqref{eq: KeyL2}, it follows that $\lims\limits_{k\rightarrow +\infty} \big(\sigma_k J_k(\beta)\big)=0$. This allows us to split the right-hand side of \eqref{eq: KeyL4} into two series, one of which is telescoping. Using that $\sigma_0=0$, we come then to
\begin{gather}
\int\limits_{B(x,r)}\Psi_{\beta}(x')\rho_{\delta}(x,x') \D \mathfrak{m}(x')\leq \sum\limits_{k=0}^{+\infty} J_k(\beta)(\sigma_{k+1}-\sigma_k).
\label{eq: KeyL5}
\end{gather}

Notice now that the functions $\Psi_{\beta}$, $\beta\in \mathbb{R}_{>0}$, increase pointwise on $B(x,r)$ as $\beta\rightarrow +\infty$ to $\chi_E(\bdot)\big(\mathsf{d}_f(x,\bdot)\big)^p$. This observation, in turn, enables us, by the monotone convergence theorem for integrals, to infer that the numbers $J_k(\beta)$, $\beta\in \mathbb{R}_{>0}$, increase to $J_k$ as $\beta\rightarrow +\infty$. Now we pass to the limit as $\beta\rightarrow +\infty$ in \eqref{eq: KeyL5} and, taking into account that the sequence $(\sigma_k)_{k\in \mathbb{N}_0}$ is non-decreasing, apply the monotone convergence theorem for series. If we juxtapose \eqref{eq: KeyL5},\eqref{eq: KeyL1}, \eqref{eq: MNotPU}, and \eqref{eq: REnD}, we end up with the upper estimate from \eqref{eq: KeyL}.

To obtain the lower estimate from \eqref{eq: KeyL} and finish the proof, one need only to provide an analogous reasoning with respect to the inequalities
\begin{gather*}
\int\limits_{B(x,r)}\Psi_{\beta}(x')\rho_{\delta}(x,x') \D \mathfrak{m}(x')\geq \sum\limits_{k=0}^{+\infty} \sigma_{k}\big(J_{k}(\beta)-J_{k+1}(\beta)\big), \quad \beta\in \mathbb{R}_{>0} ,
\end{gather*}
which take place also due to condition iii) of mollifiers.
\end{proof}
\end{lem}

The following lemma is about a metric-valued analog to the well-known fundamental pointwise inequalities that are implied by doubling and Poincar{\'e} conditions. Although they are typically formulated in terms of the maximal functions, we want to exploit a more general form, namely, the one that involves the Riesz potentials. As a consequence of these inequalities, we obtain also that maps of finite Cheeger energies have the Luzin-Lipschitz property in the sense of \cref{def: LuzLip}. For the lemma, we recall \cref{def: MetCh}, \cref{prop: LipSup}, \cref{prop: PoinIn}, \cref{prop: LebDif}, and assertion iii) from \cref{prop: RegProp}, as well as the notation given in \eqref{eq: AOp}-\eqref{eq: MOp}.
\begin{lem}
\label{lem: ChKorSchU}
Let $p\in \mathbb{R}_{\geq 1}$, suppose $\mathrm{X}$ is internally doubling and supports an internal $p$-Poincar{\'e} inequality. Then there exists a constant $c_U\in \mathbb{R}_{>0}$ depending only on $p,C_D, C_P$ such that the following holds. Let $\Omega\subseteq \mathrm{X}$ be an open set, let $(\mathrm{Y}, \mathsf{d}_{\mathrm{Y}})$ be a metric space, let $f\in \mathfrak{B}(\Omega,\mathrm{Y})$. Fix a Borel set $E\Subset \Omega$, put $\lambda\coloneqq 2\lambda_P$ and $R\coloneqq \tfrac{1}{8\lambda}\dist\big(E,\overline{\mathrm{X}}\backslash \Omega\big)$. There are $Q_f\subseteq \Omega$ with $\mathfrak{m}(\Omega\backslash Q_f)=0$ and $r_U\in (0,R)$ such that, for all $r\in (0,r_U]$, $x\in (E\cap Q_f)$, and $x'\in \big(B(E,r)\cap Q_f\big)$, we have
\begin{equation}
\big(\mathsf{d}_f(x,x')\big)^p \leq c_U r^p\Big(\mathcal{R}_p[f](x,\lambda r)+\mathcal{R}_p[f](x',\lambda r)\Big).
\label{eq: ChKorSchU0}
\end{equation}
Suppose, in addition, that $\mathrm{E}_p[f](\Omega)<+\infty$. Then $f\big|_E$ has the Luzin-Lipschitz property. 
\begin{proof}
Apply \cref{prop: LipSup}, so let $Q_{SV}\subseteq \Omega$ and $(\phi_j)_{j\in \mathbb{N}}\subseteq \mathrm{BLip}_1(\mathrm{Y})$ be as there. Fix $j\in \mathbb{N}$, put $u_j\coloneqq \phi_j\circ f$, let $Q_j$ be the set of Lebesgue points of $u_j$. Since $u_j\in \Leb^1_{Loc}(\Omega)$, we have $\mathfrak{m}(\Omega\backslash Q_j)=0$ by \cref{prop: LebDif}. Put then $Q_f\coloneqq Q_{SV}\cap \bigcap\limits_{j\in \mathbb{N}}Q_j$, whence we have $\mathfrak{m}(\Omega\backslash Q_f)=0$.

Put $r_U\coloneqq \frac{1}{64}\min\big\{R, R_D(E),R_P(E)\big\}$. From \cref{prop: LocCom}, it follows that the balls $B(x,r)$, with $x\in B(E,r_U)$ and $r\in (0,\lambda r_U]$, are compact. Since such balls are also contained in $\Omega$, we can apply \cref{prop: PoinIn} to them.

Pick $r\in (0,r_U]$, $x\in (E\cap Q_f)$, $x'\in \big(B(x,r)\cap Q_f\big)$, and $j\in \mathbb{N}$. For any $k\in \mathbb{N}_0$, we put $B_k\coloneqq B\big(x,\tfrac{r}{2^{k-1}}\big)$ and $ B'_k\coloneqq B\big(x',\tfrac{r}{2^{k-1}}\big)$. Applying the standard telescoping argument, together with the triangle, doubling, and Poincar{\'e} inequalities, we get
\begin{gather}
\mathsf{d}_{u_j}(x,x')\leq \big|\langle u_j\rangle_{B'_1}-\langle u_j\rangle_{B_0}\big|+\big|\langle u_j\rangle_{B_0}-u(x)\big|
+\big|\langle u_j\rangle_{B'_1}-u_j(x')\big|\leq \notag\\ 
\leq \fint\limits_{B_1'} \Big|u_j(x'')-\langle u_j\rangle_{B_0}\Big|\D\mathfrak{m}(x'')+\sum\limits_{k=0}^{+\infty} \Big|\langle u_j\rangle_{B_k}-\langle u_j\rangle_{B_{k+1}}\Big|+\sum\limits_{k=0}^{+\infty}\Big|\langle u_j\rangle_{B_{k+1}'}-\langle u_j\rangle_{B_{k+2}'}\Big|\leq\notag\\
\leq C_D^2\fint\limits_{B_0} \Big|u_j(x'')-\langle u_j\rangle_{B_0}\Big|\D\mathfrak{m}(x'')+\sum\limits_{k=0}^{+\infty} \fint\limits_{B_{k+1}}\Big|\langle u_j\rangle_{B_k}-u_j(x'')\Big|\D\mathfrak{m}(x'')
+\notag\\
+\sum\limits_{k=0}^{+\infty} \fint\limits_{B_{k+2}'}\Big|\langle u_j\rangle_{B_{k+1}'}-u_j(x'')\Big|\D\mathfrak{m}(x'')\leq\notag\\
\leq 2C_D^2 \sum\limits_{k=0}^{+\infty} \fint\limits_{B_k}\Big|u_j(x'')-\langle u_j\rangle_{B_k}\Big|\D\mathfrak{m}(x'') +C_D \sum\limits_{k=0}^{+\infty} \fint\limits_{B_k'}\Big|u_j(x'')-\langle u_j\rangle_{B_{k}'}\Big|\D\mathfrak{m}(x'')\leq\notag\\
\leq 4C_D^2C_P r\sum\limits_{k=0}^{+\infty} \tfrac{1}{2^k}\bigg(\Big(\mathcal{A}_p[f]\big(x,\tfrac{\lambda r}{2^k}\big)\Big)^{\frac{1}{p}}+\Big(\mathcal{A}_p[f]\big(x',\tfrac{\lambda r}{2^k}\big)\Big)^{\frac{1}{p}}\bigg),
\label{eq: PointInRealRiesz1}
\end{gather}
where we also used that $\mathrm{E}_p[f]$ majorizes $\mathrm{E}_p[u_j]$ on any set, which is by \cref{def: MetCh}.

Notice that H{\"o}lder's inequality for series gives for all sequences $(a_k)_{k\in \mathbb{N}_0}, ({a_k'})_{k\in \mathbb{N}_0}\subseteq [0,+\infty]$ that
\begin{gather*}
\Bigg(\sum\limits_{k=0}^{+\infty} \tfrac{1}{2^k} (a_k+{a_k'})\Bigg)^p=
\Bigg(\sum\limits_{k=0}^{+\infty} \big(\tfrac{2}{3}\big)^{\frac{k}{p}} (a_k+{a_k'}) \big(\tfrac{3}{2^{p+1}}\big)^{\frac{k}{p}}\Bigg)^p\leq\notag\\
\leq \Bigg(1+\sum\limits_{k=1}^{+\infty} \big(\tfrac{3}{2^{p+1}}\big)^{\frac{k}{p-1}}\Bigg)^{p-1}\sum\limits_{k=0}^{+\infty}\big(\tfrac{2}{3}\big)^k (a_k+{a_k'})^p
\leq\notag\\
\leq \biggg(\frac{2}{1-\big(\tfrac{3}{2^{p+1}}\big)^{\frac{1}{p-1}}}\biggg)^{p-1}\sum\limits_{k=0}^{+\infty}\big(\tfrac{2}{3}\big)^k \big(a_k^p+{a_k'}^p\big).
\label{eq: ChKorSchU4}
\end{gather*}
Bearing the notation from \eqref{eq: ROp} in mind, once we take, for any $k\in \mathbb{N}_0$, $\Big(\mathcal{A}_p[f]\big(x,\tfrac{\lambda r}{2^k}\big)\Big)^{\frac{1}{p}}$ as $a_k$ and $\Big(\mathcal{A}_p[f]\big(x',\tfrac{\lambda r}{2^k}\big)\Big)^{\frac{1}{p}}$ as ${a_k'}$, the above estimate, when combined with \eqref{eq: PointInRealRiesz1}, yields
\begin{equation*}
\big(\mathsf{d}_{u_j}(x,x')\big)^p\leq c_U r^p\Big(\mathcal{R}_p[f](x,\lambda r)+\mathcal{R}_p[f](x',\lambda r)\Big),
\end{equation*}
where $c_U\in \mathbb{R}_{>0}$ is some constant depending only on $p$, $C_D$, $C_P$. To derive \eqref{eq: ChKorSchU0} and finish the proof of the first part, it remains to take the supremum from the above expression over $j\in \mathbb{N}$ and use \eqref{eq: LipSup} from \cref{prop: LipSup}.

Now we move to the second part. Since $\mathrm{X}$ is separable, one can find a set $(x_m)_{m\in \mathbb{N}}\subseteq Q_f$ that is dense in $Q_f$. For each $m\in \mathbb{N}$, put
\begin{equation*}
E_m\coloneqq \bigg\{x\in \Big(E\cap B\big(x_m,\tfrac{r_U}{2}\big)\cap Q_f\Big)\Bigm| \mathcal{M}_p[f](x,\lambda r_U)\leq m\bigg\}.
\end{equation*}
Then, by the fact that $(x_m)_{m\in \mathbb{N}}$ is dense in $Q_f$ and by assertion iii) in \cref{prop: RegProp}, we have $\mathfrak{m}\bigg(E\Big\backslash \bigcup\limits_{m\in \mathbb{N}} E_m\bigg)=0$. Moreover, picking $m\in \mathbb{N}$ and $x,x'\in E_m$, applying \eqref{eq: ChKorSchU0} with $ \mathsf{d}(x,x')$ taken as $r$, and using the definition of $E_m$, one can easily verify that
\begin{equation*}
\mathsf{d}_f(x,x')\leq (2mc_U)^{\frac{1}{p}}\mathsf{d}(x,x').
\end{equation*}
This observation shows that $f\big|_{E_m}$ is Lipschitz for any $m\in \mathbb{N}$. Thus, $f\big|_E$ has the Luzin-Lipschitz property. The proof is now complete.
\end{proof}
\end{lem}

The lemma below demonstrates how to construct for a given function a member of an approximating sequence in the sense of \cref{def: ChEnerg}. Its, rather standard, proof is based on a use of a discrete convolution constructed via Lipschitz partition of unity as in \cref{prop: RegProp}.
\begin{lem}\label{lem: ChKorSchL}
Let $p\in \mathbb{R}_{\geq 1}$, suppose $\mathrm{X}$ is internally doubling. Then there exists a constant $c_L\in \mathbb{R}_{>0}$ depending only on $p,C_D$ such that the following holds. Let $\Omega\subseteq \mathrm{X}$ be an open set, let $u\in \Leb^1_{Loc}(\Omega)$. Fix an open set $\Omega'\Subset \Omega$. There exists $r_L\in \mathbb{R}_{>0}$ such that, for every $r\in (0,r_L]$, one can find $u^r\in \mathrm{Lip}_{Loc}(\Omega')$ fulfilling the conditions below:
\begin{gather}
\big|u^r(x)-u(x)\big|^p\leq c_L r^p\mathsf{ks}_p[u,\Omega](x,r) \quad \text{for $\mathfrak{m}$-a.e. $x\in \Omega'$};\label{eq: ChKorSchL0}\\
 \big(\mathrm{lip}[u^r](x)\big)^p\leq c_L\mathsf{ks}_p[u,\Omega](x,r)\quad \text{for $\mathfrak{m}$-a.e. $x\in \Omega'$}.
\label{eq: ChKorSchL00}
\end{gather}
\begin{proof}
Taking $\Omega'$ as $S$, we apply \cref{prop: PartUnit}, so let $r_{PU}$ be as there. Put then $r_L\coloneqq\min\Big\{\dist\big(\Omega', \overline{\mathrm{X}}\backslash \Omega\big), R_D(\Omega'),r_{PU}\Big\}$ and fix $r\in (0,r_L]$. Considering $\tfrac{r}{8}$ as $r$ for the proposition, we find families $(x_i)_{i\in \mathbb{N}}$ and $(\varphi_i)_{i\in \mathbb{N}}$ with the corresponding properties, which we refer to below as properties i)-iii). 

Define a function $u^r\colon \Omega' \to \mathbb{R}$ by
\begin{equation*}
u^r(x)\coloneqq\sum\limits_{i\in \mathbb{N}} \varphi_i(x)\langle u\rangle_{B(x_i,\frac{r}{8})}.
\end{equation*}
This function is well posed by the following reasons. For any $i\in \mathbb{N}$, the ball $B(x_i,r)$ is compact, which is due to \cref{prop: LocCom}, whence $u\in \Leb^1\big(B(x_i,r)\big)$. In turn, properties i) and ii) imply that the set of all $i\in \mathbb{N}$ with $\varphi_i(x)\neq 0$ is finite for any $x\in \Omega'$.

Using properties i)-iii) and Jensen's and the doubling inequalities, one can, for any $x\in \Omega'$, get
\begin{gather*}
\big|u^r(x)-u(x)\big|^p=\Bigg|\sum\limits_{i\in \mathbb{N}}\varphi_i(x) \Big(\langle u\rangle_{B(x_i,\frac{r}{8})}-u(x)\Big)\Bigg|^p \leq\notag\\
 \leq c_{PU} \sup\limits_{i\in \mathbb{N}}\Bigg(\chi_{B(x_i, \frac{r}{8})}(x)\Big|\langle u\rangle_{B(x_i,\frac{r}{8})}-u(x)\Big|^p\Bigg)\leq\notag\\
\leq  c_{PU} C_D^2\fint\limits_{B(x,\frac{r}{4})} \big(\mathsf{d}_u(x,x')\big)^p\D\mathfrak{m}(x')\leq c_{PU} C_D^4 r^p\mathsf{ks}_p[u,\Omega](x,r),
\label{eq: BlaBla1}
\end{gather*}
which corresponds to \eqref{eq: ChKorSchL0}. Next, fix $x_0\in \Omega'$ and find $r'\in \big(0,\tfrac{r}{32}\big)$ such that $B(x_0,r')\subseteq \Omega'$. Fix $x,x'\in B(x_0,r')$. By properties ii) and iii), one can find $i_0\in \mathbb{N}$ such that $x_0\in B\big(x_{i_0}, \tfrac{r}{8}\big)$. Using then properties i)-iii), again together with the doubling inequality, we get
\begin{gather*}
\big|u^r(x')-u^r(x)\big|
=\Bigg|\sum\limits_{i\in \mathbb{N}}\big(\varphi_i(x')-\varphi_i(x)\big)\Big(\langle u\rangle_{B(x_i,\frac{r}{8})}-\langle u\rangle_{B(x_{i_0},\frac{r}{8})}\Big)\Bigg|\leq\notag\\
\leq \frac{c_{PU}^2}{r}\mathsf{d}(x,x')\sup\limits_{i\in \mathbb{N}} \bigg(\chi_{B(x_i,\frac{r}{4})}(x)\Big|\langle u\rangle_{B(x_i,\frac{r}{8})}-\langle u\rangle_{B(x_{i_0},\frac{r}{8})}\Big|\bigg)\leq \notag\\
\leq \frac{2c_{PU}^2C_D^3}{r}\mathsf{d}(x,x')\fint\limits_{B(x_0,\frac{r}{2})} \mathsf{d}_u(x_0,x'')\D\mathfrak{m}(x'')\leq  \frac{2c_{PU}^2C_D^4}{r}\mathsf{d}(x,x') \fint\limits_{B(x_0,r)} \mathsf{d}_u(x_0,x'')\D\mathfrak{m}(x'').
\end{gather*}
This inequality, together with the fact that $u\in \Leb^1\big(B(x_0,r)\big)$, which, in turn, is again due to \cref{prop: LocCom}, and the arbitrariness of $x$ and $x'$, shows that $u^r\big|_{B(x_0,r')}$ is Lipschitz. Thus, we conclude that $u^r\in \mathrm{Lip}_{Loc}(\Omega)$. Moreover, after letting $x=x_0$, dividing the expression by $\mathsf{d}(x,x')$, passing to the upper limit as $x'\to x_0$, and applying H{\"o}lder's inequality, one can easily verify that
\begin{equation*}
\big(\mathrm{lip}[u^r](x_0)\big)^p\leq (2c_{PU}^2 C_D^4)^p\mathsf{ks}_p[u,\Omega](x_0,r),
\label{eq: BlaBla2}
\end{equation*}
which corresponds to \eqref{eq: ChKorSchL00}.

Everything we established means that the conclusion of the lemma will follow once we put $c_L\coloneqq \max\big\{c_{PU}C_D^4,(2c_{PU}^2C^4_D)^p\big\}$.
\end{proof}
\end{lem}
\subsection{The first result}\label{ss: Res1} Our intention here is to establish the validity of \cref{theo: Res1}.

From this point until the end of the section, the following is implied:
\begin{itemize}
\item[$\bullet $] we are given a parameter $p\in \mathbb{R}_{\geq 1}$;
\item[$\bullet $] we assume that $\mathrm{X}$ is internally doubling;
\item[$\bullet $] we are given a $p$-admissible family $(\rho_{\delta})_{\delta\in (0,1)}$.
\end{itemize}

The following statement provides a local version of the upper estimate in \eqref{eq: Res1} from \cref{theo: Res1}. And it is worth noticing that the corresponding relation is established with no additional restrictions neither on the open set nor on the map. This yields a lower bound of a general nature for the considered energy. The proof of the theorem is based on the combination of \cref{lem: ChKorSchU} with the upper bound from \cref{lem: KeyL}. We recall assertion i) from \cref{prop: RegProp} and condition v) of mollifiers as well.
\begin{theo}\label{lem: ResUp}
Suppose $\mathrm{X}$ supports an internal $p$-Poincar{\'e} inequality. There exists a constant $C_U\in \mathbb{R}_{>0}$ depending only on $p,C_D,C_P,\lambda_P,C_M$ such the following holds. Let $\Omega\subseteq \mathrm{X}$ be an open set, let $(\mathrm{Y},\mathsf{d}_{\mathrm{Y}})$ be a metric space, let $f\in \mathfrak{B}(\Omega,\mathrm{Y})$. We have
\begin{equation}
\sup\limits_{\Omega'}\lim\limits_{r\searrow 0}\lims\limits_{\delta\searrow 0} \int\limits_{\Omega'}\int\limits_{B(x,r)} \big(\mathsf{d}_f(x,x')\big)^p\rho_{\delta}(x,x') \D \mathfrak{m}(x')\D \mathfrak{m}(x)\leq C_U \mathrm{E}_p[f](\Omega),
\label{eq: ResUp0}
\end{equation}
where the supremum is taken over all open sets $\Omega'
\Subset \Omega$.
\begin{proof}
If $\mathrm{E}_p[f](\Omega)=+\infty$, then \eqref{eq: ResUp0} holds trivially, so further we consider the converse case.

Fix an open set $\Omega'\Subset \Omega$, put $\lambda\coloneqq 4\lambda_P$. Apply \cref{lem: ChKorSchU}, with $\Omega'$ considered as $E$, so let $Q_f$ and $r_U$ be as there. Since \eqref{eq: ChKorSchU0} holds for all $x\in (\Omega'\cap Q_f)$ and $x'\in \big(B(\Omega',r_U)\cap Q_f\big)$ and since $\mathfrak{m}(\Omega\backslash Q_f)=0$, we can pick any $r\in (0,r_U)$ and integrate \eqref{eq: ChKorSchU0} with respect to $x'\in B(x,r)$ first and then with respect to $x\in \Omega'$. After that, exploiting assertion i) from \cref{prop: RegProp} several times and using the doubling inequality, we can readily check that
\begin{equation}
\int\limits_{\Omega'}\mathsf{ks}_p[f,\Omega](x,r)\D \mathfrak{m}(x)\leq 2C_D^{\log_2(8\lambda)}c_U \mathrm{E}_p[f]\big(B(\Omega',\lambda r)\big)\leq 2C_D^{\log_2(8\lambda)}c_U \mathrm{E}_p[f](\Omega).
\label{eq: ResUp1}
\end{equation}
Next, apply \cref{lem: KeyL}, taking $\Omega$ as $E$. Picking $\delta\in (0,1)$ and $r\in (0,r_U)$, letting $h=2$, and integrating the second inequality from \eqref{eq: KeyL0} with respect to $x\in \Omega'$, we come to
\begin{gather*}
\int\limits_{\Omega'}\int\limits_{B(x,r)} \Big(\mathsf{d}_f(x,x')\big)^p\rho_{\delta}(x,x') \D \mathfrak{m}(x')\D \mathfrak{m}(x)\leq \notag\\
\leq\sum\limits_{k=0}^{+\infty} \int\limits_{\Omega'} \Pi_{U,k}(\delta,x,r,2)\mathsf{ks}_p[f,\Omega]\big(x,\tfrac{r}{2^k}\big)\D \mathfrak{m}(x)\leq \notag\\
\leq \Bigg(\sum\limits_{k=0}^{+\infty}\esup\limits_{x\in \Omega'}\Pi_{U,k}(\delta,x,r,2)\Bigg)\sup\limits_{r'\in (0,r]}\int\limits_{\Omega'}\mathsf{ks}_p[f,\Omega](x,r')\D \mathfrak{m}(x).
\end{gather*}
In the above inequality, we pass to the upper limits as $\delta\searrow 0$ and then as $r\searrow 0$. Now, with respect to the expression thus obtained, we implement the upper estimate in \eqref{eq: Cond5} from condition v) of mollifiers. Combine then the result with \eqref{eq: ResUp1}, after passing there to the upper limit as $r\searrow 0$, and put $C_U\coloneqq 2C_D^{\log_2(16\lambda)}c_UC_M$. From the arbitrariness of $\Omega'$, the assertion now follows.
\end{proof}
\end{theo}

The next statement concerns a localized form of the lower estimate in \eqref{eq: Res1} from \cref{theo: Res1}. And as with the previous theorem, we do not make any assumptions on the map and on the open set, which enables us to obtain a pretty generic upper bound for the corresponding energy. The proof of the theorem relies on \cref{lem: ChKorSchL} and on the lower estimate from \cref{lem: KeyL}. We recall also \cref{def: ChEnerg} and \cref{def: MetCh}, as well as \cref{prop: SupM}.
\begin{theo}\label{lem: ResD}
There exists a constant $C_L\in \mathbb{R}_{>0}$ depending only on $p,C_D,C_M$ such the following holds. Let $\Omega\subseteq \mathrm{X}$ be an open set, let $(\mathrm{Y}, \mathsf{d}_{\mathrm{Y}}\big)$ be a metric space, let $f\in \mathfrak{B}(\Omega,\mathrm{Y})$. We have
\begin{equation}
C_L \mathrm{E}_p[f](\Omega)\leq \sup\limits_{\Omega'}\lim\limits_{r\searrow 0}\limi\limits_{\delta\searrow 0} \int\limits_{\Omega'}\int\limits_{B(x,r)} \big(\mathsf{d}_f(x,x')\big)^p\rho_{\delta}(x,x') \D \mathfrak{m}(x')\D \mathfrak{m}(x),
\label{eq: ResD0}
\end{equation}
where the supremum is taken over all open sets $\Omega'\Subset \Omega$.
\begin{proof}
Fix an open set $\Omega'\Subset \Omega$. We start by applying \cref{lem: KeyL}, so that, after putting $R\coloneqq \tfrac{1}{2}\min\Big\{\dist\big(\Omega',\overline{\mathrm{X}}\backslash \Omega\big), R_D(\Omega')\Big\}$, picking $\delta\in (0,1)$ and $r\in (0,R)$, letting $h=2$, and integrating the first inequality from \eqref{eq: KeyL0} with respect to $x\in \Omega'$, one can get
\begin{gather*}
 \int\limits_{\Omega'}\int\limits_{B(x,r)} \big(\mathsf{d}_f(x,x')\big)^p\rho_{\delta}(x,x') \D \mathfrak{m}(x')\D \mathfrak{m}(x)\geq \notag\\
\geq\sum\limits_{k=0}^{+\infty} \int\limits_{\Omega'} \Pi_{L,k}(\delta,x,r,2)\mathsf{ks}_p[f,\Omega]\big(x,\tfrac{r}{2^k}\big)\D \mathfrak{m}(x)\geq \notag\\
\geq \sum\limits_{k=0}^{+\infty}\einf\limits_{x\in \Omega'}\Pi_{L,k}(\delta,x,r,2)\int\limits_{\Omega'}\mathsf{ks}_p[f,\Omega]\big(x,\tfrac{r}{2^k}\big)\D \mathfrak{m}(x)\geq\\
\geq \Bigg(\sum\limits_{k=0}^{+\infty}\einf\limits_{x\in \Omega'}\Pi_{L,k}(\delta,x,r,2)\Bigg)\inf\limits_{r'\in (0,r]}\int\limits_{\Omega'}\mathsf{ks}_p[f,\Omega](x,r')\D \mathfrak{m}(x).
\label{eq: ResD1}
\end{gather*}
Now we pass here to the lower limits as $\delta\searrow 0$ and then as $r\searrow 0$. After this, using the lower estimate in \eqref{eq: Cond5} from condition v) of mollifiers, we obtain that
\begin{gather}
\limi\limits_{r\searrow 0} \int\limits_{\Omega'}\mathsf{ks}_p[f,\Omega](x,r)\D \mathfrak{m}(x)\leq C_M\lim\limits_{r\searrow 0}\limi\limits_{\delta\searrow 0}\int\limits_{\Omega'}\int\limits_{ B(x,r)} \big(\mathsf{d}_f(x,x')\big)^p\rho_{\delta}(x,x') \D \mathfrak{m}(x')\D \mathfrak{m}(x).
\label{eq: ResD2}
\end{gather}

Taking into account assertion ii) from \cref{prop: SupM}, consider now an arbitrary pairwise disjoint family $(\Omega_j)_{j\in \mathbb{N}}$ of open subsets of $\Omega'$ and an arbitrary family $(\phi_j)_{j\in \mathbb{N}}\subseteq \mathrm{BLip}_1(\mathrm{Y})$, put then $u_j\coloneqq \phi_j\circ f$ for any $j\in \mathbb{N}$. Estimating from below the integral over $\Omega'$ in the left-hand side of \eqref{eq: ResD2} in terms of the series of the same integrals but over the sets $\Omega_j$, $j\in \mathbb{N}$, applying then Fatou's lemma for series, and using that the functions $\phi_j$, $j\in \mathbb{N}$, are $1$-Lipschitz, we deduce that
\begin{gather}
\sum\limits_{j=1}^{+\infty} \limi\limits_{r\searrow 0} \int\limits_{\Omega_j}\mathsf{ks}_p[u_j,\Omega](x,r)\D \mathfrak{m}(x)\leq  \notag\\
\leq C_M\lim\limits_{r\searrow 0}\limi\limits_{\delta\searrow 0}\int\limits_{\Omega'}\int\limits_{ B(x,r)} \big(\mathsf{d}_f(x,x')\big)^p\rho_{\delta}(x,x') \D \mathfrak{m}(x')\D \mathfrak{m}(x).
\label{eq: ResD3}
\end{gather}
If the right-hand side of \eqref{eq: ResD3} is infinite, then \eqref{eq: ResD0} is satisfied in a trivial way, so we consider further the opposite case. Consequently, we have
\begin{equation}
\limi\limits_{r\searrow 0}\int\limits_{\Omega_j}\mathsf{ks}_p[u_j,\Omega](x,r)\D \mathfrak{m}(x)<+\infty\quad\text{for any $j\in \mathbb{N}$}.
\label{eq: ResD4}
\end{equation}
Now we fix $j\in \mathbb{N}$ and, taking $\Omega_j$ as $\Omega'$ and $u_j$ as $u$, we apply \cref{lem: ChKorSchL}, so let $r_L$ be as there and let $u^r$, for any $r\in (0,r_L)$, denote the corresponding function constructed with respect to $r$. Since \eqref{eq: ResD4} holds, there is a sequence $(r_n)_{n\in \mathbb{N}}\subseteq (0,r_L)$ converging to zero as $n\rightarrow +\infty$ such that
\begin{gather}
\limi\limits_{r\searrow 0}\int\limits_{\Omega_j}\mathsf{ks}_p[u_j,\Omega](x,r)\D \mathfrak{m}(x)=\lim\limits_{n\rightarrow +\infty}\int\limits_{\Omega_j}\mathsf{ks}_p[u_j,\Omega](x,r_n)\D \mathfrak{m}(x)<+\infty.
\label{eq: ResD5}
\end{gather}
By the very formulation of \cref{lem: ChKorSchL}, we have that $u^{r_n}\in \mathrm{Lip}_{Loc}(\Omega_j)$. Then, taking \eqref{eq: ResD5} into account, we infer via \eqref{eq: ChKorSchL0} that $(u^{r_n})_{n\in \mathbb{N}}$ converges to $u$ in $\Leb^p(\Omega_j)$, and hence in $\Leb^p_{Loc}(\Omega_j)$, as $n\rightarrow +\infty$. Now we juxtapose \eqref{eq: ChEnerg} from \cref{def: ChEnerg} with \eqref{eq: ChKorSchL00} from \cref{lem: ChKorSchL} and conclude, using also \eqref{eq: ResD5}, that
\begin{equation*}
\mathrm{E}_p[u_j](\Omega_j)\leq c_L\limi\limits_{r\searrow 0} \int\limits_{\Omega_j}\mathsf{ks}_p[u_j,\Omega](x,r)\D \mathfrak{m}(x).
\end{equation*}

Since $j$ is arbitrary, it follows from \eqref{eq: ResD3} and the above inequality that
\begin{equation*}
\sum\limits_{j=1}^{+\infty}\mathrm{E}_p[u_j](\Omega_j)\leq c_LC_M\lim\limits_{r\searrow 0}\limi\limits_{\delta\searrow 0}\int\limits_{\Omega'}\int\limits_{B(x,r)} \big(\mathsf{d}_f(x,x')\big)^p\rho_{\delta}(x,x') \D \mathfrak{m}(x')\D \mathfrak{m}(x).
\end{equation*}
This, together with from assertion ii) from \cref{prop: SupM} and the arbitrariness of the families $(\Omega_j)_{j\in \mathbb{N}}$ and $(\phi_j)_{j\in \mathbb{N}}$, leads to
\begin{equation*}
\mathrm{E}_p[f](\Omega')\leq c_LC_M\lim\limits_{r\searrow 0}\limi\limits_{\delta\searrow 0}\int\limits_{\Omega'}\int\limits_{B(x,r)}\big(\mathsf{d}_f(x,x')\big)^p\rho_{\delta}(x,x')\D \mathfrak{m}(x')\D \mathfrak{m}(x).
\end{equation*}
Now, after taking the supremum over all open sets $\Omega'\Subset \Omega$, applying assertion i) in \cref{prop: SupM}, and putting $C_L\coloneqq \frac{1}{c_LC_M}$, we get \eqref{eq: ResD0}, which finishes the proof.
\end{proof}
\end{theo}

Now we are ready to finish the proof of \cref{theo: Res1}. Suppose $\mathrm{X}$ supports an internal $p$-Poincar{\'e} inequality. Let $\Omega\Subset \mathrm{X}$ be an open set having the strong $p$-extension property in the sense of \cref{def: ExrProp}, let $(\mathrm{Y}, \mathsf{d}_{\mathrm{Y}})$ be a metric space, let $f\in \Leb^p(\Omega, \mathrm{Y})$. Fix first an arbitrary metric space $(\mathrm{Z}, \mathsf{d}_{\mathrm{Z}})$ and an arbitrary $\mathrm{Z}$-extension $F$ of $f$ belonging to $\mathfrak{B}(\mathrm{X}, \mathrm{Z})$. Pick any $R\in \mathbb{R}_{>0}$ and fix any open set $\Omega_0\subseteq B(\Omega,R)$ with $\Omega\Subset \Omega_0$. Considering $\Omega_0$ as $\Omega$, $\mathrm{Z}$ as $\mathrm{Y}$, and $F\big|_{\Omega_0}$ as $f$, we apply \cref{lem: ResUp}. After this, we use \cref{lem: ResD}, taking $\mathrm{Z}$ as $\mathrm{Y}$ and $F\big|_{\Omega}$ as $f$. Combining all this and putting $C\coloneqq \max\big\{C_U, \tfrac{1}{C_L}\big\}$, one can derive that
\begin{gather}
\frac{1}{C} \mathrm{E}_p[F](\Omega)\leq \lim\limits_{r\searrow 0}\limi\limits_{\delta\searrow 0} \int\limits_{\Omega}\int\limits_{\Omega\cap B(x,r)} \big(\mathsf{d}_F(x,x')\big)^p\rho_{\delta}(x,x') \D \mathfrak{m}(x')\D \mathfrak{m}(x)\leq \notag\\
\leq \lim\limits_{r\searrow 0} \lims\limits_{\delta\searrow 0} \int\limits_{\Omega}\int\limits_{\Omega\cap B(x,r)}\big(\mathsf{d}_F(x,x')\big)^p\rho_{\delta}(x,x') \D \mathfrak{m}(x')\D \mathfrak{m}(x)\leq C \mathrm{E}_p[F]\big(B(\Omega,R)\big) .
\label{eq: PartRes0}
\end{gather}
Since $F$ is an extension of $f$, we can replace $F$ by $f$ in the first three terms above. This trivially leads to \eqref{eq: Res1} in the case $\mathrm{E}_p[f](\Omega)=+\infty$. Suppose now that $\mathrm{E}_p[f](\Omega)<+\infty$. Then, since $\Omega$ has the strong $p$-extension property, we may, without loss of generality, require $F$ to fulfill \eqref{eq: ZerEn}. Hence, once we pass to the limit as $R\searrow 0$ in \eqref{eq: PartRes0}, its right-hand side, up to the constant there, becomes equal to $\mathrm{E}_p[F](\Omega)$, which, in turn, coincides with $\mathrm{E}_p[f](\Omega)$. For \eqref{eq: Res1} to be established, it remains now to use that $f\in \Leb^p(\Omega, \mathrm{Y})$ and apply \eqref{eq: BExt0} from \cref{lem: BExt}. Thus, \cref{theo: Res1} is completely proved.
\subsection{The second result}\label{ss: Res2} Within this section, we provide the proof for \cref{theo: Res2}.

The following suppositions are common for the entire section:
\begin{itemize}
\item[$\bullet $] we are given a parameter $p\in \mathbb{R}_{\geq 1}$;
\item[$\bullet $] we assume that $\mathrm{X}$ is internally doubling;
\item[$\bullet $] we assume that $\mathrm{X}$ supports an internal $p$-Poincar{\'e} inequality.
\end{itemize}

It is in the following theorem, where we basically establish the relations, or, more precisely, their localized variants, from \eqref{eq: Res2} and \eqref{eq: Res2'}. But in it, instead of imposing conditions on the underlying space, we require a map to satisfy a property mimicking the one guaranteed by \cref{theo: ApprMetrDiff}. All this is worth pointing out since such a property can be established in more general situations, rather than only when the source space is strongly rectifiable. In addition, as an intermediate but independently valuable result, we show that under the above-mentioned property there exists, for every $p\in \mathbb{R}_{\geq 1}$, the $p$-energy density, in the sense of \cref{def: EnD}, for any well-behaved map. When combined with \cref{theo: ApprMetrDiff}, this fact actually becomes an extension of the corresponding result from \cite{GT21} to the case $p=1$. But once again, we believe that this statement is important in isolation from strongly rectifiable spaces as well. For the proof of theorem, we recall \cref{prop: GenLebTh}, \cref{def: RelGr}, \cref{prop: LebDif}, assertions i) and ii) from \cref{prop: RegProp}, the conditions of mollifiers from \cref{def: Mol}, \cref{lem: ChKorSchU}, \cref{lem: KeyL}, as well as the notation given in \eqref{eq: AOp}, \eqref{eq: ROp}, and \eqref{eq: REnD}.
\begin{theo}\label{lem: EnDEx}
Let $\Omega\subseteq \mathrm{X}$ be an open set, let $(\mathrm{Y},\mathsf{d}_{\mathrm{Y}})$ be a metric space, let $f\in \Leb^p_{Loc}(\Omega,\mathrm{Y})\cap\mathrm{S}^{1,p}(\Omega,\mathrm{Y})$. Suppose there exists a function $e\in \mathfrak{B}(\Omega,\mathbb{R}_{\geq 0})$ with the property that, for $\mathfrak{m}$-a.e. $x\in \Omega$, one can find a sequence $(V_{x,n})_{n\in \mathbb{N}}$ of Borel subsets of $ \Omega$ having $x$ as a density point such that
\begin{equation}
\lims\limits_{n\rightarrow +\infty}\lims\limits_{r\searrow 0} \Big| \mathsf{ks}_p[f,V_{x,n}](x,r)-e(x)\Big|=0.
\label{eq: EndEx000}
\end{equation}
Then the following assertions hold:
\begin{itemize}
\item[$\rm i)$] the function $e$ is the $p$-energy density of $f$;
\item[$\rm ii)$] given a strongly $p$-admissible family $(\rho_{\delta})_{\delta\in (0,1)}$ and an open set $\Omega'\Subset\Omega$, we have
\begin{gather}
\lims\limits_{r \searrow 0}\lims\limits_{\delta \searrow 0}\biggg|\int\limits_{\Omega'\cap B(x,r)}\big(\mathsf{d}_f(x,x')\big)^p\rho_{\delta}(x,x')\D\mathfrak{m}(x')-\Theta_M(x)\mathsf{e}_p[f](x)\biggg|=0\quad \text{for $\mathfrak{m}$-a.e. $x\in \Omega'$},\label{eq: ConvEnD0}\\
\lims\limits_{r\searrow 0}\lims\limits_{\delta \searrow 0}\int\limits_{\Omega'}\biggg|\int\limits_{\Omega'\cap B(x,r)}\big(\mathsf{d}_f(x,x')\big)^p\rho_{\delta}(x,x')\D\mathfrak{m}(x')-\Theta_M(x)\mathsf{e}_p[f](x)\biggg|\D\mathfrak{m}(x)=0.\label{eq: ConvEnD00}
\end{gather}
\end{itemize}
\begin{proof}
We start with assertion i).

Let $Q_0\subseteq \Omega $ be the set of all points satisfying the condition from the theorem hypothesis, then $\mathfrak{m}(\Omega\backslash Q_{0})=0$. Let $r_U$ and $Q_f$ be as in \cref{lem: ChKorSchU}, with $\Omega'$ considered as $E$, put $\lambda\coloneqq 2\lambda_P$ and $R\coloneqq \tfrac{1}{8\lambda}\min\big\{R_D(\Omega'),r_U\big\}$. As belonging to $\mathrm{S}^{1,p}(\Omega,\mathrm{Y})$, the map $f$ has the minimal $p$-weak gradient, so we can put $G\coloneqq |\nabla f|_p^p$. Let $Q_G$ be the set of all Lebesgue points of $G$. Since $G\in \Leb^1(\Omega)$, from \cref{prop: LebDif} it follows that $\mathfrak{m}(\Omega\backslash Q_G)=0$. Put $Q\coloneqq  (Q_{0}\cap Q_f\cap Q_G)$, whence $\mathfrak{m}(\Omega\backslash Q)=0$.

We show first that the functions $\mathsf{ks}_p[f,\Omega](\bdot,r)$, $r\in (0,R)$, converge $\mathfrak{m}$-a.e. on $\Omega'$ to $e$ as $r\searrow 0$. Fix $x\in Q$, $n\in \mathbb{N}$, and $r\in (0,R)$. By exploiting \eqref{eq: ChKorSchU0} and applying the triangle inequality multiple times, we can obtain that
\begin{gather}
\Big|\mathsf{ks}_p[f, V_{x,n}](x,r)-\mathsf{ks}_p[f, \Omega](x,r)\Big|\leq \notag\\
\leq  \frac{1}{\mathfrak{m}\big(B(x,r)\big)}\int\limits_{B(x,r)\backslash V_{x,n}}\frac{\big(\mathsf{d}_f(x,x')\big)^p}{r^p}\D\mathfrak{m}(x')+\frac{\mathfrak{m}\big(B(x,r)\backslash V_{x,n}\big)}{\mathfrak{m}\big(B(x,r)\cap V_{x,n}\big)}\mathsf{ks}_p[f,\Omega](x,r)\leq \notag\\ \leq \frac{2c_U\mathfrak{m}\big(B(x,r)\backslash V_{x,n}\big)}{\mathfrak{m}\big(B(x,r)\cap V_{x,n}\big)}\mathcal{R}_p[f](x,\lambda r) +\notag\\
+\frac{c_U}{\mathfrak{m}\big(B(x,r)\big)}\int\limits_{B(x,r)\backslash V_{x,n}} \mathcal{R}_p[f](x',\lambda r)\D \mathfrak{m}(x')+\notag\\
+'\frac{c_U\mathfrak{m}\big(B(x,r)\backslash V_{x,n}\big)}{\mathfrak{m}\big(B(x,r)\cap V_{x,n}\big)}\fint\limits_{B(x,r)} \mathcal{R}_p[f](x',\lambda r)\D \mathfrak{m}(x')\leq \notag\\
\leq \frac{2c_U\mathfrak{m}\big(B(x,r)\backslash V_{x,n}\big)}{\mathfrak{m}\big(B(x,r)\cap V_{x,n}\big)}\Big|\mathcal{R}_p[f](x,\lambda r)-G(x)\Big|+\frac{2c_U\mathfrak{m}\big(B(x,r)\backslash V_{x,n}\big)}{\mathfrak{m}\big(B(x,r)\cap V_{x,n}\big)}G(x)+\notag\\
+c_U \fint\limits_{B(x,r)} \Big|\mathcal{R}_p[f](x',\lambda r)-G(x)\Big|\D \mathfrak{m}(x')+\frac{c_U\mathfrak{m}\big(B(x,r)\backslash V_{x,n}\big)}{\mathfrak{m}\big(B(x,r)\big)}G(x)+\notag\\
+\frac{c_U\mathfrak{m}\big(B(x,r)\backslash V_{x,n}\big)}{\mathfrak{m}\big(B(x,r)\cap V_{x,n}\big)}\fint\limits_{B(x,r)} \Big|\mathcal{R}_p[f](x',\lambda r)-G(x)\Big|\D \mathfrak{m}(x')+\frac{c_U\mathfrak{m}\big(B(x,r)\backslash V_{x,n}\big)}{\mathfrak{m}\big(B(x,r)\cap V_{x,n}\big)}G(x).
\label{eq: EnDEx1}
\end{gather}
We claim now that all the terms in the right-hand side of the above expression tend to zero as $r \searrow 0$. We argue as follows to show this. Since $x$ is a density point of $V_{x,n}$, we have
\begin{equation}
\lims\limits_{r\searrow 0} \frac{\mathfrak{m}\big(B(x,r)\backslash V_{x,n}\big)}{\mathfrak{m}\big(B(x,r)\big)}=\lims\limits_{r\searrow 0}\frac{\mathfrak{m}\big(B(x,r)\backslash V_{x,n}\big)}{\mathfrak{m}\big(B(x,r)\cap V_{x,n}\big)}=0.
\label{eq: EnDEx2}
\end{equation}
Next, keeping in mind the notation given in \eqref{eq: AOp} and \eqref{eq: ROp} and the facts that $x$ is a Lebesgue point of $G$ and that $G$ is the Radon-Nikodym derivative of $\mathrm{E}_p[f]$ with respect to $\mathfrak{m}\resmes_{\Omega}$, we get
\begin{gather}
\lims\limits_{r \searrow 0}\Big|\mathcal{R}_p[f](x,\lambda r)-G(x)\Big|\leq \tfrac{1}{3}\lims\limits_{r \searrow 0}\sum\limits_{k=0}^{+\infty} \big(\tfrac{2}{3}\big)^k\fint\limits_{B\big(x',\frac{\lambda r}{2^k}\big)} \mathsf{d}_G(x,x')\D \mathfrak{m}(x')\leq \notag\\
\leq \lims\limits_{r \searrow 0}\sup\limits_{r'\in (0,\lambda r]}\fint\limits_{B(x, r')}\mathsf{d}_G(x,x')\D \mathfrak{m}(x')=0.
\label{eq: EnDEx3} 
\end{gather}
Besides this, the same reasoning as above, together with assertion i) from \cref{prop: RegProp} and the doubling inequality, ensures that
\begin{gather}
\lims\limits_{r \searrow 0}\fint\limits_{B(x,r)}\Big|\mathcal{R}_p[f](x',\lambda r)-G(x)\Big|\D \mathfrak{m}(x')\leq \notag\\
\leq \tfrac{1}{3}\lims\limits_{r \searrow 0}\sum\limits_{k=0}^{+\infty} \big(\tfrac{2}{3}\big)^k\fint\limits_{B(x,r)}\fint\limits_{B\big(x',\frac{\lambda r}{2^k}\big)} \mathsf{d}_G(x,x'')\D \mathfrak{m}(x'')\D \mathfrak{m}(x) \leq\notag\\
\leq \lims\limits_{r \searrow 0}\sum\limits_{k=0}^{+\infty} \big(\tfrac{2}{3}\big)^k\frac{1}{\mathfrak{m}\big(B(x,r)\big)}\int\limits_{B(x,2\lambda r)} \mathsf{d}_G(x,x')\D \mathfrak{m}(x')\leq \notag\\
\leq C_D^{\log_2 (4\lambda)}\lims\limits_{r \searrow 0} \sup\limits_{r'\in (0,2\lambda r]}\fint\limits_{B(x,r')}\mathsf{d}_G(x,x')\D \mathfrak{m}(x')=0.
\label{eq: EnDEx4}
\end{gather}
Now we combine \eqref{eq: EnDEx2}-\eqref{eq: EnDEx4} with \eqref{eq: EnDEx1} and conclude that the claim is valid. By \eqref{eq: EndEx000}, it then follows via the arbitrariness of $x$ that 
\begin{equation}
\lims\limits_{r\searrow 0}\Big|\mathsf{ks}_p[f, \Omega](x,r)-e(x)\Big|=0 \quad \text{for $\mathfrak{m}$-a.e. $x\in \Omega'$},
\label{eq: EnDEx5}
\end{equation}
as desired.

Now we want to verify that the functions $\mathsf{ks}_p[f,\Omega](\bdot,r)$, $r\in (0,R)$, converge to $e$ in $\Leb^1(\Omega')$ as $r\searrow 0$. Picking $x\in Q$ and $r\in (0,R)$, we integrate \eqref{eq: PointInRealRiesz1} with respect to $x'\in B(x,r)$ and, applying assertion i) from \cref{prop: RegProp} together with the doubling inequality as previously, get
\begin{equation*}
\mathsf{ks}_p[f,\Omega](x,r)\leq c_U\mathcal{R}_p[f](x,\lambda r)+c_UC_D^{\log_2(4\lambda)}\mathcal{A}_p[f](x,2\lambda r).
\end{equation*}
Since \eqref{eq: EnDEx5} holds and since the right-sided expression above converges in $\Leb^1(\Omega')$ by assertion ii) from \cref{prop: RegProp}, we are exactly in position to apply \cref{prop: GenLebTh}, whence the stated convergence follows.

Since $\Omega'$ is arbitrary, the conclusion of assertion i) follows.

We now turn to assertion ii).

Let $R$ be as above with respect to $\Omega'$. Pick $\delta\in (0,1)$, $r\in(0,R)$, and $h\in \mathbb{R}_{>1}$. Define functions $g_{L,h,r,\delta},g_{U,h,r,\delta}\in \mathfrak{B}(\Omega',\mathbb{R}_{\geq 0})$ as
\begin{gather*}
g_{L,h,r,\delta}\coloneqq \biggg|\Theta_M(x) \mathsf{e}_p[f](x)-\sum\limits_{k=0}^{+\infty}\Pi_{L,k}(\delta,x,r,h)\mathsf{ks}_p[f,\Omega']\big(x,\tfrac{r}{h^k}\big)\biggg|,\\
g_{U,h,r,\delta}\coloneqq \biggg|\Theta_M(x) \mathsf{e}_p[f](x)-\sum\limits_{k=0}^{+\infty}\Pi_{U,k}(\delta,x,r,h)\mathsf{ks}_p[f,\Omega']\big(x,\tfrac{r}{h^k}\big)\biggg|,
\end{gather*}
respectively. We easily deduce from \eqref{eq: KeyL0} that
\begin{gather}
\biggg|\int\limits_{\Omega'\cap B(x,r)}\big(\mathsf{d}_f(x,x')\big)^p\rho_{\delta}(x,x')\D\mathfrak{m}(x')-\Theta_M(x)\mathsf{e}_p[f](x)\biggg|\leq g_{L,h,r,\delta}(x)+g_{U,h,r,\delta}(x)
\label{eq: ConvEnD10}
\end{gather}
for any $x\in \Omega'$. Fixing $x\in \Omega'$ and applying some simple transformations to $g_{U,h,r,\delta}(x)$, we obtain that
\begin{gather}
g_{U,h,r,\delta}(x)\leq \mathsf{e}_p[f](x)\Bigg|\Theta_M(x)-\sum\limits_{k=0}^{+\infty}\Pi_{U,k}(\delta,x,r,h)\Bigg|+\notag\\
+\sum\limits_{k=0}^{+\infty}\Pi_{U,k}(\delta,x,r,h)\Big| \mathsf{ks}_p[f,\Omega']\big(x,\tfrac{r}{h^k}\big)-\mathsf{e}_p[f](x)\Big|\leq\notag\\
\leq \mathsf{e}_p[f](x)\Bigg|\Theta_M(x)-\sum\limits_{k=0}^{+\infty}\Pi_{U,k}(\delta,x,r,h)\Bigg|+\notag\\
+\sup\limits_{r'\in (0,r]}\Big| \mathsf{ks}_p[f,\Omega'](x,r')-\mathsf{e}_p[f](x)\Big|\sum\limits_{k=0}^{+\infty}\sup\limits_{x'\in \Omega'}\Pi_{U,k}(\delta,x',r,h).
\label{eq: ConvEnD2}
\end{gather}
Now we pass here to the upper limits as $\delta\searrow 0$, then as $r\searrow 0$, and finally as $h\searrow 1$ and conclude, via condition vi) of mollifiers and \cref{def: EnD}, that
\begin{equation*}
\lims\limits_{h\searrow 1}\lims\limits_{r\searrow 0}\lims\limits_{\delta\searrow 0}g_{U,h,r,\delta}(x)=0\quad \text{for $\mathfrak{m}$-a.e. $x\in \Omega'$}.
\label{eq: ConvEnD6}
\end{equation*}
After carrying out similar reasoning, but with respect to $g_{L,h,r,\delta}$, we also have 
\begin{equation*}
\lims\limits_{h\searrow 1}\lims\limits_{r\searrow 0}\lims\limits_{\delta\searrow 0}g_{L,h,r,\delta}(x)=0\quad \text{for $\mathfrak{m}$-a.e. $x\in \Omega'$}.
\label{eq: ConvEnD7}
\end{equation*}
Combining these facts with \eqref{eq: ConvEnD10}, we come to \eqref{eq: ConvEnD0}, so the first part is verified.

Next, picking $h\in \mathbb{R}_{>1}$, $r\in (0,R)$, and $\delta\in (0,1)$, integrating $g_{U,h,r,\delta}$ over $\Omega'$, and following the same strategy as in \eqref{eq: ConvEnD2}, we get
\begin{gather}
\int\limits_{\Omega'}g_{U,h,r,\delta}(x)\D \mathfrak{m}(x) \leq \int\limits_{\Omega'}\mathsf{e}_p[f](x)\Bigg|\Theta_M(x)-\sum\limits_{k=0}^{+\infty}\Pi_{U,k}(\delta,x,r,h)\Bigg|\D \mathfrak{m}(x)+\notag\\
+\biggg(\sup\limits_{r'\in (0,r]}\int\limits_{\Omega'}\Big| \mathsf{ks}_p[f,\Omega'](x,r')-\mathsf{e}_p[f](x)\Big|\D \mathfrak{m}(x)\biggg)\sum\limits_{k=0}^{+\infty}\sup\limits_{x'\in \Omega'}\Pi_{U,k}(\delta,x',r,h).
\label{eq: ConvEnD20}
\end{gather}
We claim that both terms in the right-hand side of the above expression converge to zero if we pass to the upper limits as $\delta\searrow 0$, then as $r \searrow 0$, and finally as $h\searrow 0$. To show this, we act as follows. Begin with the first term. By condition vi) of mollifiers and \cref{def: EnD}, we are able to interchange, via revers Fatou's lemma, the limits as $\delta\searrow 0$ and as $r\searrow 0$ with the integral sign. After that, we use conditions vi) and vii) of mollifiers, \cref{def: EnD}, and the dominated convergence theorem in order to perform such a manipulation once again, but with the limit as $h\searrow 1$. Thereby, the resulted limit of the first term is indeed equal to zero.
Continue with the second term. Notice that
\begin{equation*}
\mathsf{ks}_p[f,\Omega'](x,r)\leq \mathsf{ks}_p[f,\Omega](x,r)\quad \text{for all $x\in \Omega'$ and $r\in (0,R)$}.
\end{equation*}
By \cref{def: EnD}, the right-sided expression above converges to $\mathsf{e}_p[f]$ in $\Leb^1(\Omega')$ as $r\searrow 0$, while the left-sided expression clearly converges $\mathfrak{m}$-a.e. on $\Omega'$ as $r\searrow 0$ to $\mathsf{e}_p[f]$. In light of these observations, we are again in position to apply \cref{prop: GenLebTh} and conclude that the functions $\mathsf{ks}_p[f,\Omega'](\bdot,r)$, $r\in (0,R)$, converge to $\mathsf{e}_p[f]$ in $\Leb^1(\Omega')$ as $r\searrow 0$. From this fact and condition vi) of mollifiers, it yields that the limit of the second term also equals zero. Consequently, we have
\begin{equation*}
\lims\limits_{h\searrow 1}\lims\limits_{r\searrow 0}\lims\limits_{\delta\searrow 0}\int\limits_{\Omega'}g_{U,h,r,\delta}(x)\D\mathfrak{m}(x)=0.
\end{equation*}
And as previously, one can verify in an analogous way that
\begin{equation*}
\lims\limits_{h\searrow 1}\lims\limits_{r\searrow 0}\lims\limits_{\delta\searrow 0}\int\limits_{\Omega'}g_{L,h,r,\delta}(x)\D\mathfrak{m}(x)=0.
\end{equation*}
It remains only to combine both of these facts with \eqref{eq: ConvEnD10} and obtain \eqref{eq: ConvEnD00}, which is the second part of the assertion.

Thus, we complete the proof of assertion ii), and hence the theorem. 
\end{proof}
\end{theo}

We are able now to establish \cref{theo: Res2}. Suppose $\mathrm{X}$ is strongly rectifiable. Let $(\rho_{\delta})_{\delta\in (0,1)}$ be a strongly $p$-admissible family. Let $\Omega\Subset \mathrm{X}$ be an open set, let $(\mathrm{Y}, \mathsf{d}_{\mathrm{Y}})$ be a metric space, let $f\in \mathrm{W}^{1,p}(\mathrm{X},\mathrm{Y})$. Taking $\mathrm{X}$ as $\Omega$ and $\Omega$ as $E$, we apply \cref{lem: ChKorSchU} and conclude that $f\big|_{\Omega}$ has the Luzin-Lipschitz property, which enables us to use \cref{theo: ApprMetrDiff}. After combining it with assertion ii) from \cref{lem: EnDEx}, with $\mathrm{X}$ and $\Omega$ considered as $\Omega$ and $\Omega'$, respectively, putting $\epsilon_p[f]\coloneqq \Theta_M\mathsf{e}_p[f]$, and applying \eqref{eq: BExt0} and \eqref{eq: BExt00} from \cref{lem: BExt}, we end up with \eqref{eq: Res2} and \eqref{eq: Res2'}, so the proof of \cref{theo: Res2} is now fully complete.
\section{Corollaries}\label{ss: App}
\subsection{Examples of mollifiers}\label{ss: ExMol}
Now we show that \cref{theo: Res1} and \cref{theo: Res2} are not vacuously true due to the non-existence of the corresponding families of mollifiers. 

There are several acceptable examples of mollifiers that are typically considered in the literature. And we are about to verify that our results are perfectly applicable at least to those of them that are discussed in \cite{LPZ22}. So let us list them first. Let $p\in \mathbb{R}_{\geq 1}$. Define, for any $\delta\in (0,1)$, functions $\rho^0_{p,\delta},\rho^1_{p,\delta},\rho^2_{p,\delta},\rho^3_{p,\delta}\in \mathfrak{B}\big(\mathrm{X}\times \mathrm{X},\mathbb{R}_{\geq 0}\big)$ by
\begin{gather}
\rho^0_{p, \delta}(x,x')\coloneqq\frac{\delta}{\big(\mathsf{d}(x,x')\big)^{p(1-\delta)} \mathfrak{m}\Big(B\big(x,4\mathsf{d}(x,x')\big)\Big)},\label{eq: Mol0}\\
\rho^1_{p,\delta}(x,x')\coloneqq\frac{\chi_{B(x,\delta)}(x')}{{\delta}^p \mathfrak{m}\big(B(x,\delta)\big)},\label{eq: Mol1}\\
\rho^2_{p,\delta}(x,x')\coloneqq\frac{\chi_{B(x,\delta)}(x')}{\big(\mathsf{d}(x,x')\big)^p \mathfrak{m}\big(B(x,\delta)\big)},\label{eq: Mol2}\\
\rho^3_{p,\delta}(x,x')\coloneqq\frac{\chi_{B(x,\delta)}(x')}{{\delta}^p \mathfrak{m}\Big(B\big(x,\mathsf{d}(x,x')\big)\Big)}\label{eq: Mol3},
\end{gather}
respectively. Note that the family $(\rho^3_{p,\delta})_{\delta\in (0,1)}$, while not being considered in the mentioned source, is provided here as fitting into a "pattern" of how mollifiers should look like and as also suitable for our context.

\begin{rema}
To clarify for the reader, the reason we added the, rather obscure, number four in \eqref{eq: Mol0} is that without it, even provided that the underlying space is internally doubling, the respective family may not fulfill condition i) of mollifiers as decreasing too slowly at infinity. In general, one can put there any number greater than two. But once the space is globally doubling, or at least satisfies some limitation on the growth rate of the measures of balls at large scales, the number can be made equal to one. And in such a case, all the subsequent statements regarding these mollifiers, but with the number four eliminated, are valid as well. This circumstance can be derived from the corresponding proofs.
\end{rema}

Now we are ready to provide some transparent sufficient conditions for the families defined via \eqref{eq: Mol0}-\eqref{eq: Mol3} to be admissible and strongly admissible on the given space.
\begin{lem}\label{prop: AmdCond} Let $p\in \mathbb{R}_{\geq 1}$. Suppose $\mathrm{X}$ is internally doubling. Then each family from $(\rho^0_{p, \delta})_{\delta\in (0,1)}$, $(\rho^1_{p,\delta})_{\delta\in (0,1)}$, $(\rho^2_{p, \delta})_{\delta\in (0,1)}$, $(\rho^3_{p, \delta})_{\delta\in (0,1)}$ is $p$-admissible, and, moreover, the related constant $C_M$ can be chosen so that it depends only on $p$ and $C_D$. Suppose, in addition, that there exists a function $D\in \mathfrak{B}(\mathrm{X},\mathbb{R}_{>0})$ such that for $\mathfrak{m}$-a.e. $x\in \mathrm{X}$ one has
\begin{equation}
\lim\limits_{r\searrow 0}\frac{\mathfrak{m}\big(B(x,hr)\big)}{\mathfrak{m}\big(B(x,r)\big)}=h^{D(x)}\quad \text{for any $h\in \mathbb{R}_{>0}$}.
\label{eq: StrAdmC}
\end{equation}
Then each of the listed families is strongly $p$-admissible, and, moreover, the related function $\Theta_M$ is equal to the corresponding one among the functions $\Theta_0,\Theta_1, \Theta_2,\Theta_3\in \mathfrak{B}(\mathrm{X},\mathbb{R}_{\geq 0})$ given by
\begin{gather}
\Theta_0(x)\coloneqq \frac{D(x)+p}{4^{D(x)}p},\quad \Theta_1(x)\coloneqq 1, \quad\Theta_2(x)\coloneqq \frac{D(x)+p}{D(x)},\quad\Theta_3(x)\coloneqq \frac{D(x)+p}{p},
\label{eq: StrAdmC0}
\end{gather}
respectively.
\begin{proof}
The proofs for all given families, though pretty laborious, are basically analogous and consist in a careful check of all the conditions from \cref{def: Mol}. In this regard, for the sake of brevity, we provide details for the family $(\rho^0_{p,\delta})_{\delta\in (0,1)}$ only. For any $\delta\in (0,1)$, put $\rho^0_{p,\delta}\coloneqq \rho_{\delta}$.

Fix $x\in \mathrm{X}$, $r\in \mathbb{R}_{>0}$, and $\delta\in (0,1)$. Put $A_k\coloneqq B\big(x,{3}^{k+1}r\big)\backslash B\big(x,{3}^{k}r\big)$ for any $k\in \mathbb{N}_0$ and $K\coloneqq \sup\big\{k\in \mathbb{N}\mid \mathfrak{m}(A_k)<+\infty\big\}$. It can be checked then that
\begin{gather*}
\int\limits_{\mathrm{X}\backslash B(x,r)}\big(\rho_{\delta}(x,x')+\rho_{\delta}(x,x')\big)\D \mathfrak{m}(x')\leq \notag\\ \leq 2\delta\sum\limits_{k=0}^{K}\int\limits_{A_k}\frac{1}{\big(\mathsf{d}(x,x')\big)^{p(1-\delta)}\mathfrak{m}\Big(B\big(x,{3}\mathsf{d}(x,x')\big)\Big)}\D \mathfrak{m}(x')\leq \notag\\
\leq 2\delta\sum\limits_{k=0}^{K}\frac{\mathfrak{m}(A_k)}{({3}^k r)^{p(1-\delta)} \mathfrak{m}\big(B(x,3^{k+1}r)\big)}\leq \frac{2\delta {3}^{p(1-\delta)}}{r^{p(1-\delta)}\big({3}^{p(1-\delta)}-1\big)}.
\end{gather*}
We now take here the supremum over $x\in \mathrm{X}$, pass to the upper limits as $\delta\searrow 0$ and then as $r\searrow 0$, which gives \eqref{eq: Cond1}, so condition i) of mollifiers is satisfied.

Next, we clearly have
\begin{equation}
\sigma(\delta,x,r)\leq \frac{\delta}{r^{p(1-\delta)}\mathfrak{m}\big(B(x,4 r)\big)}\quad\text{for all $\delta\in (0,1)$, $x\in \mathrm{X}$, and $r\in \mathbb{R}_{>0}$},
\label{eq: ExMol1}
\end{equation}
whence, after passing to the upper limits as $\delta\searrow 0$ and then as $r\searrow 0$, we come to \eqref{eq: Cond2}, so condition ii) of mollifiers is satisfied as well.

From the very definition of the considered family, it easily follows that condition iii) of mollifiers is valid for it.

After multiplying the inequality from \eqref{eq: ExMol1} by $\Big(\mathfrak{m}\big(B(x,r)\big)r^p\Big)$ and passing in the resulted expression to the upper limit as $r\searrow 0$, we come to \eqref{eq: Cond4}, whence we deduce that condition iv) of mollifiers holds too.

Fix a nonempty set $S\Subset\mathrm{X}$. Put $R\coloneqq \tfrac{1}{2}\min\Big\{\dist\big(S, \overline{\mathrm{X}}\backslash \mathrm{X}\big), R_D(S)\Big\}$. Then, since we assumed that our space is connected and of positive diameter, it can be seen that $\big(B(x,r')\backslash B(x,r)\big)\neq \emptyset$ for all $x\in S$ and $r,r'\in (0,R)$ with $r<r'$. Fix $\delta\in (0,1)$, $x\in S$, $r\in (0,R)$, and $h\in \mathbb{R}_{>1}$. Simple calculations show that
\begin{gather}
\Pi_{L,k}(\delta,x,r,h)=\begin{dcases}\tfrac{\delta r^{p\delta}}{h^{2p(1-\delta)}(h^{p\delta})^k}\frac{\mathfrak{m}\Big(B\big(x,\tfrac{r}{h^{k+1}}\big)\Big)}{\mathfrak{m}\Big(B\big(x,\tfrac{4 r}{h^{k-1}}\big)\Big)}\biggg(\frac{h^{p(1-\delta)}\mathfrak{m}\Big(B\big(x,\tfrac{4 r}{h^{k-1}}\big)\Big)}{\mathfrak{m}\Big(B\big(x,\tfrac{4 r}{h^{k}}\big)\Big)}-1\biggg),\quad k\neq 0,\\
\tfrac{\delta r^{p\delta}}{h^{p(1-\delta)}}\frac{\mathfrak{m}\Big(B\big(x,\tfrac{r}{h}\big)\Big)}{\mathfrak{m}\big(B(x,4 r\big)},\quad \text{otherwise},
\end{dcases}\label{eq: ExMol2}\\
\Pi_{U,k}(\delta,x,r,h)=\begin{dcases}\tfrac{\delta  r^{p\delta}}{h^{p(1-\delta)}(h^{p\delta})^k}\frac{\mathfrak{m}\Big(B\big(x,\tfrac{r}{h^{k}}\big)\Big)}{\mathfrak{m}\Big(B\big(x,\tfrac{4 r}{h^{k-1}}\big)\Big)}\biggg(\frac{h^{p(1-\delta)}\mathfrak{m}\Big(B\big(x,\tfrac{4 r}{h^{k-1}}\big)\Big)}{\mathfrak{m}\Big(B\big(x,\tfrac{4 r}{h^{k}}\big)\Big)}-1\biggg),\quad k\neq 0,\\
\delta r^{p\delta}\frac{\mathfrak{m}\big(B(x,r )\big)}{\mathfrak{m}\big(B(x,4r )\big)},\quad \text{otherwise}.
\end{dcases}\label{eq: ExMol3}
\end{gather}
From the doubling inequality, it follows for any $k\in \mathbb{N}_0$ that
\begin{gather*}
\frac{h^{p(1-\delta)}-1}{C_D^{\log_2(8 h^2)}}\tfrac{\delta}{h^{2p(1-\delta)}}\big(\tfrac{r}{h^{k}}\big)^{p\delta}\leq \Pi_{L,k}(\delta,x,r,h)\leq \Pi_{U,k}(\delta,x,r,h)\leq \delta \big(\tfrac{r}{h^{k}}\big)^{p\delta} h^{p(1-\delta)}C_D^{\log_2(2h)},
\end{gather*}
whence we have
\begin{gather}
\frac{h^{p(1-\delta)}-1}{C_D^{\log_2(8 h^2)}}\frac{\delta r^{p\delta}}{h^{2p(1-\delta)}}\frac{1}{1-\frac{1}{h^{p\delta}}}\leq \sum\limits_{k=0}^{+\infty} \inf\limits_{x\in \Omega} \Pi_{L,k}(\delta,x,r,h)\leq \notag\\
\leq  \sum\limits_{k=0}^{+\infty} \sup\limits_{x\in \Omega} \Pi_{U,k}(\delta,x,r,h)\leq \frac{\delta r^{p\delta}h^{p(1-\delta)}C_D^{\log_2(2h)}}{1-\frac{1}{h^{p\delta}}}.
\label{eq: ExMol4}
\end{gather}
After picking $h= 2$, passing in the above expression to the upper limits as $\delta \searrow 0$ and then as $r\searrow 0$ and putting $C_M\coloneqq \frac{p2^{2p}C_D^{4}}{\ln(2)}$, we conclude that condition v) of mollifiers holds. So the first part of the statement is verified.

We move to the second part. Put $\Theta\coloneqq \Theta_0$.

Notice that \eqref{eq: ExMol4} allows us to deduce immediately that \eqref{eq: Cond6} is met, whence the validity of condition vi) of mollifiers follows.

Fix $x\in \mathrm{X}$ and $\delta\in (0,1)$, find $r\in \mathbb{R}_{>0}$ with $\mathfrak{m}\big(B(x,r)\big)<+\infty$, and $h\in \mathbb{R}_{>1}$. Keeping \eqref{eq: ExMol3} and \eqref{eq: StrAdmC0} in mind, it follows via the triangle inequality that
\begin{gather*}
\biggg|\Theta(x)-\sum\limits_{k=1}^{+\infty} \Pi_{U,k}(\delta,x,r,h)\biggg|\leq \Theta(x)\bigg|1-\frac{\delta r^{p\delta}h^{p\delta}p\ln(h)}{h^{p(1-\delta)}(h^{p\delta}-1)}\bigg|+\notag\\
+\delta (h^{p\delta}-1)\sum\limits_{k=1}^{+\infty} \big(\tfrac{r}{h^k}\big)^{p\delta}\Biggg|\frac{\mathfrak{m}\Big(B\big(x,\tfrac{r}{h^{k}}\big)\Big)}{\mathfrak{m}\Big(B\big(x,\tfrac{4 r}{h^{k}}\big)\Big)}\Biggg|+\notag\\
+ \sum\limits_{k=1}^{+\infty}\frac{\delta p\ln(h)}{h^{p(1-\delta)}}\big(\tfrac{r}{h^k}\big)^{p\delta}\Biggg|\frac{1}{p\ln(h)}\frac{\mathfrak{m}\Big(B\big(x,\tfrac{r}{h^{k}}\big)\Big)}{\mathfrak{m}\Big(B\big(x,\tfrac{4 r}{h^{k-1}}\big)\Big)}\biggg(\frac{h^{p}\mathfrak{m}\Big(B\big(x,\tfrac{4 r}{h^{k-1}}\big)\Big)}{\mathfrak{m}\Big(B\big(x,\tfrac{4 r}{h^{k}}\big)\Big)}-1\biggg)-\Theta(x)\Biggg|.
\end{gather*}
Pass in the above inequality to the upper limits as $\delta\searrow 0$ and then as $r\searrow 0$. After performing all necessary estimates, applying \eqref{eq: StrAdmC}, and taking into account that $\lims\limits_{r\searrow 0}\lims\limits_{\delta\searrow 0} \Pi_{U,0}(\delta,x,r,h)=0$, it follows for $\mathfrak{m}$-a.e. $x\in \mathrm{X}$ that
\begin{gather}
\lims\limits_{r\searrow 0}\lims\limits_{\delta\searrow 0}\biggg|\Theta(x)-\sum\limits_{k=0}^{+\infty} \Pi_{U,k}(\delta,x,r,h)\biggg|\leq \notag\\
\leq \Theta(x)\big|1-\tfrac{1}{h^{p}}\big|+\tfrac{1}{h^p}\lims\limits_{r\searrow 0}\biggg|\frac{1}{p\ln(h)}\frac{\mathfrak{m}\big(B(x,r)\big)}{\mathfrak{m}\big(B(x,4 h r)\big)}\Bigg(\frac{h^{p}\mathfrak{m}\big(B(x,4 h r)\big)}{\mathfrak{m}\big(B(x, 4 r)\big)}-1\Bigg)-\Theta(x)\biggg|=\notag\\
=\Theta(x)\big|1-\tfrac{1}{h^{p}}\big|+\tfrac{1}{h^p}\Big|\tfrac{1}{p\ln(h)(4 h)^{D(x)}}(h^{p+D(x)}-1)-\Theta(x)\Big|=\frac{D(x)+p}{p4^{D(x)}}\big|1-\tfrac{1}{h^{p}}\big|+\notag\\
+\tfrac{1}{ph^p (4 h)^{D(x)}}\Big|\tfrac{1}{\ln(h) }(h^{p+D(x)}-1)-D(x)-p\Big|+\tfrac{h^{D(x)}-1}{ph^p(4 h)^{D(x)}}\big|p+D(x)\big|.
\label{eq: ExMol5}
\end{gather}
Notice that the doubling inequality, together with \eqref{eq: StrAdmC}, trivially implies that $D(x)\leq \log_2(C_D)$ for $\mathfrak{m}$-a.e. $x\in \mathrm{X}$. Therefore, after taking the essential supremum from \eqref{eq: ExMol5} over $x\in \mathrm{X}$, passing to the upper limit as $h\searrow 1$, and applying some standard reasoning from calculus, we end up with
\begin{equation*}
\lims\limits_{h\searrow 1}\esup\limits_{x\in \mathrm{X}}\lims\limits_{r\searrow 0}\lims\limits_{\delta\searrow 0}\biggg|\Theta(x)-\sum\limits_{k=0}^{+\infty} \Pi_{U,k}(\delta,x,r,h)\biggg|=0.
\end{equation*}
Following the same ideas but with the use of the quantities from \eqref{eq: ExMol2}, we get
\begin{equation*}
\lims\limits_{h\searrow 1}\esup\limits_{x\in \mathrm{X}}\lims\limits_{r\searrow 0}\lims\limits_{\delta\searrow 0}\biggg|\Theta(x)-\sum\limits_{k=0}^{+\infty} \Pi_{L,k}(\delta,x,r,h)\biggg|=0.
\end{equation*}
All this ensures the fulfillment of \eqref{eq: Cond7}, and hence of condition vii) of mollifiers.

Thus, the proof of the second part, as well as the proposition itself, is now complete.
\end{proof}
\end{lem}

Summing up, by the lemma just proved, the families defined via \eqref{eq: Mol0}-\eqref{eq: Mol3} are indeed admissible on internally doubling spaces, and hence are appropriate candidates to be considered for \cref{theo: Res1}. For the same to be the case with respect to \cref{theo: Res2}, it only remains to notice that by \eqref{eq: HausDens0} from \cref{theo: HausDens}, one can take the function $\mathrm{Dim}$ as $D$ to fulfill \eqref{eq: StrAdmC} whenever the source space is strongly rectifiable.
\subsection{Final relations} We conclude our paper by formulating explicitly the outcome of the application of our results to the discussed mollifiers.

Combining \cref{theo: Res1} with \cref{prop: AmdCond}, we receive the following corollary about an equivalence of Cheeger energies to several various quantities obtained as the limits of nonlocal functionals.
\begin{cor}
Let $p\in \mathbb{R}_{\geq 1}$, suppose $\mathrm{X}$ is internally doubling and supports an internal $p$-Poincar{\'e} inequality.  There exist constants $C_1,C_2,C_3,C_4,C_5\in \mathbb{R}_{>0}$ with $C_1\leq C_2\leq C_3\leq C_4\leq C_5$ depending only on $p$, $C_D$, $C_P$, $\lambda_P$ such that the following holds. Let $\Omega\Subset\mathrm{X}$ be an open set having the strong $p$-extension property, let $(\mathrm{Y},\mathsf{d}_{\mathrm{Y}})$ be a metric space, let $f\in \Leb^p(\Omega,\mathrm{Y})$. The estimates below hold:
\begin{gather}
\mathrm{E}_p[f](\Omega)\leq C_1\limi\limits_{s\nearrow 1}\int\limits_{\Omega\times \Omega} \frac{(1-s)\big(\mathsf{d}_f(x,x')\big)^p}{\big(\mathsf{d}(x,x')\big)^{ps}\mathfrak{m}\Big(B\big(x,4\mathsf{d}(x,x')\big)\Big)}\D (\mathfrak{m}\otimes \mathfrak{m})(x,x')\leq \notag\\
\leq C_1\lims\limits_{s\nearrow 1}\int\limits_{\Omega\times \Omega} \frac{(1-s)\big(\mathsf{d}_f(x,x')\big)^p}{\big(\mathsf{d}(x,x')\big)^{ps}\mathfrak{m}\Big(B\big(x,4\mathsf{d}(x,x')\big)\Big)}\D (\mathfrak{m}\otimes \mathfrak{m})(x,x')\leq \notag \\
\leq C_2\limi\limits_{r\searrow 0}\int\limits_{\Omega}\int\limits_{\Omega\cap B(x,r)} \frac{\big(\mathsf{d}_f(x,x')\big)^p}{r^p\mathfrak{m}\big(B(x,r)\big)}\D\mathfrak{m}(x')\D\mathfrak{m}(x)\leq \notag\\
\leq C_2\lims\limits_{r\searrow 0}\int\limits_{\Omega}\int\limits_{\Omega\cap B(x,r)} \frac{\big(\mathsf{d}_f(x,x')\big)^p}{r^p\mathfrak{m}\big(B(x,r)\big)}\D\mathfrak{m}(x')\D\mathfrak{m}(x)\leq \notag\\
\leq C_3\limi\limits_{r\searrow 0}\int\limits_{\Omega}\int\limits_{\Omega\cap B(x,r)} \frac{\big(\mathsf{d}_f(x,x')\big)^p}{\big(\mathsf{d}(x,x')\big)^p\mathfrak{m}\big(B(x,r)\big)}\D\mathfrak{m}(x')\D\mathfrak{m}(x)\leq \notag\\
\leq C_3\lims\limits_{r\searrow 0}\int\limits_{\Omega}\int\limits_{\Omega\cap B(x,r)} \frac{\big(\mathsf{d}_f(x,x')\big)^p}{\big(\mathsf{d}(x,x')\big)^p\mathfrak{m}\big(B(x,r)\big)}\D\mathfrak{m}(x')\D\mathfrak{m}(x)\leq \notag\\
\leq C_4\limi\limits_{r\searrow 0}\int\limits_{\Omega}\int\limits_{\Omega\cap B(x,r)} \frac{\big(\mathsf{d}_f(x,x')\big)^p}{r^p\mathfrak{m}\Big(B\big(x,\mathsf{d}(x,x')\big)\Big)}\D\mathfrak{m}(x')\D\mathfrak{m}(x)\leq \notag
\\
\leq C_4\lims\limits_{r\searrow 0}\int\limits_{\Omega}\int\limits_{\Omega\cap B(x,r)} \frac{\big(\mathsf{d}_f(x,x')\big)^p}{r^p\mathfrak{m}\Big(B\big(x,\mathsf{d}(x,x')\big)\Big)}\D\mathfrak{m}(x')\D\mathfrak{m}(x)\leq C_5\mathrm{E}_p[f](\Omega).
\label{eq: FEq}
\end{gather}
\end{cor}

The combination of \cref{theo: Res2}, \cref{prop: AmdCond}, \cref{theo: ApprMetrDiff}, and \cref{theo: HausDens}, yields the corollary below, within which several exact formulas concerning the limits of nonlocal functionals are listed.
\begin{cor}
Suppose $\mathrm{X}$ is strongly rectifiable. Let $p\in \mathbb{R}_{\geq 1}$, suppose $\mathrm{X}$ is internally doubling and supports an internal $p$-Poincar{\'e} inequality. Let $\Omega\Subset \mathrm{X}$ be an open set, let $(\mathrm{Y},\mathsf{d}_{\mathrm{Y}})$ be a metric space, let $f\in \mathrm{W}^{1,p}(\mathrm{X},\mathrm{Y})$. The following equalities hold:
\begin{gather}
\lim\limits_{s\nearrow 1}\int\limits_{\Omega\times \Omega} \frac{(1-s)\big(\mathsf{d}_f(x,x')\big)^p}{\big(\mathsf{d}(x,x')\big)^{ps}\mathfrak{m}\Big(B\big(x,4\mathsf{d}(x,x')\big)\Big)}\D (\mathfrak{m}\otimes \mathfrak{m})(x,x')=\notag\\
=\int\limits_{\Omega}\frac{\mathrm{Dim}(x)+p}{4^{\mathrm{Dim}(x)}p} \fint\limits_{B^{\mathrm{Dim}(x)}}\big(\mathsf{md}_x[f](x')\big)^p\D \mathcal{L}^{\mathrm{Dim}(x)}(x')\D \mathfrak{m}(x);\label{eq: EF0}\\
\lim\limits_{r\searrow 0}\int\limits_{\Omega}\int\limits_{\Omega\cap B(x,r)} \frac{\big(\mathsf{d}_f(x,x')\big)^p}{r^p\mathfrak{m}\big(B(x,r)\big)}\D\mathfrak{m}(x')\D\mathfrak{m}(x)=\notag\\
=\int\limits_{\Omega} \fint\limits_{B^{\mathrm{Dim}(x)}}\big(\mathsf{md}_x[f](x')\big)^p\D \mathcal{L}^{\mathrm{Dim}(x)}(x')\D \mathfrak{m}(x);\label{eq: EF1}\\
\lim\limits_{r\searrow 0}\int\limits_{\Omega}\int\limits_{\Omega\cap B(x,r)} \frac{\big(\mathsf{d}_f(x,x')\big)^p}{\big(\mathsf{d}(x,x')\big)^p\mathfrak{m}\big(B(x,r)\big)}\D\mathfrak{m}(x')\D\mathfrak{m}(x)=\notag\\
=\int\limits_{\Omega}\frac{\mathrm{Dim}(x)+p}{\mathrm{Dim}(x)} \fint\limits_{B^{\mathrm{Dim}(x)}}\big(\mathsf{md}_x[f](x')\big)^p\D \mathcal{L}^{\mathrm{Dim}(x)}(x')\D \mathfrak{m}(x);\label{eq: EF2}\\
\lim\limits_{r\searrow 0}\int\limits_{\Omega}\int\limits_{\Omega\cap B(x,r)} \frac{\big(\mathsf{d}_f(x,x')\big)^p}{r^p\mathfrak{m}\Big(B\big(x,\mathsf{d}(x,x')\big)\Big)}\D\mathfrak{m}(x')\D\mathfrak{m}(x)=\notag\\
=\int\limits_{\Omega}\frac{\mathrm{Dim}(x)+p}{p} \fint\limits_{B^{\mathrm{Dim}(x)}}\big(\mathsf{md}_x[f](x')\big)^p\D \mathcal{L}^{\mathrm{Dim}(x)}(x')\D \mathfrak{m}(x) \label{eq: EF3}.
\end{gather}
\end{cor}
\begin{rema}
Let us provide now some conclusive comments on the following circumstance. Unlike the rough BBM-type estimates in \eqref{eq: FEq}, the precise BBM-type formulas in \eqref{eq: EF0}-\eqref{eq: EF3} involve, instead of the Cheeger energies, objects that are expressed in terms of the metric differentials. At the same time, the right-hand side, up to the constant there, of the original BBM formula in \eqref{eq: Seminorm} is known to be exactly the Cheeger energy in its essence. So it makes sense to ask whether the quantities of both types are proportional in a universal way, i.e. with a coefficient independent of maps. The established formulas themselves point to an obvious obstacle to this, namely the dimensional inconstancy. And this is indeed the case as an example proposed in \cite{G22} illustrates. But it turns out that even if we limit ourselves to spaces of constant dimension, the positive answer to the question can be guaranteed in the real-valued case only. For details on this phenomenon, we send the reader to \cite[Section 4]{GT21}. Therefore, it is the appearance of the right-sided expressions in the formulas presented in \eqref{eq: EF0}-\eqref{eq: EF3} that should be considered as the most appropriate one within our context.
\end{rema}
\printbibliography
\end{document}